\documentclass{elsart}
\usepackage{graphicx,psfrag}
\usepackage{amsmath}
\usepackage{amssymb}
\usepackage{calligra}
\begin{document}
\begin{frontmatter}

\title{Efficient implicit integration for finite-strain viscoplasticity with a nested multiplicative split}
\author{A.V. Shutov\corauthref{cor}},
\corauth[cor]{Corresponding author. Tel.: +7-383-333-1746;
fax: +7-383-333-1612.}
\ead{alexey.v.shutov@gmail.com}

\address{Lavrentyev Institute of Hydrodynamics\\
Pr. Lavrentyeva 15, 630090 Novosibirsk, Russia, \\
Novosibirsk State University\\
Pirogova 2, 630090 Novosibirsk, Russia}

\begin{abstract}

An efficient and reliable stress computation algorithm is presented, which
is based on implicit integration of the local evolution
equations of multiplicative finite-strain plasticity/viscoplasticity.
The algorithm is illustrated by an example involving
a combined nonlinear isotropic/kinematic hardening;
numerous backstress tensors are employed
for a better description of the material behavior.

The considered material model
exhibits the so-called weak invariance
under arbitrary isochoric changes of the reference configuration, and
the presented algorithm retains this useful property.
Even more: the weak invariance serves
as a guide in constructing this algorithm.
The constraint of inelastic incompressibility is exactly preserved as well.
The proposed method is first-order accurate.

Concerning the accuracy of the stress computation, the new algorithm
is comparable to the Euler Backward method with a subsequent correction of incompressibility (EBMSC)
and the classical exponential method (EM).
Regarding the computational efficiency, the new algorithm
is superior to the EBMSC and EM.
Some accuracy tests are presented using parameters of the aluminum alloy 5754-O and the 42CrMo4 steel.
FEM solutions of two boundary value problems using MSC.MARC are presented to show
the correctness of the numerical implementation.

\end{abstract}
\end{frontmatter}


\section{Introduction}

In many practical applications, metal forming simulations are implemented to
estimate the induced plastic anisotropy, the residual stresses and the magnitude of spring back.
It is known that the numerical
computations accounting for these effects should employ constitutive models of finite strain plasticity/viscoplasticity with
combined nonlinear isotropic/kinematic hardening \cite{LaurentGreeze2009, LaurentGreeze2010, Firat2012}.
One of such models was proposed by Shutov and Krei\ss ig within
the phenomenological context of multiplicative plasticity \cite{ShutovKrVisc}.
An important constitutive assumption of this model is the nested multiplicative split of the deformation gradient:
Within the first split, which is a classical one (cf. \cite{Bilby, Kroener}), the deformation
gradient is multiplicatively decomposed into the elastic and inelastic parts.
Within the further splits, which were introduced by Lion in \cite{LionIJP} to account for a nonlinear kinematic hardening,
the inelastic part is decomposed into
conservative (energetic) and dissipative parts.

Although there is a big number of approaches to
finite strain elasto-plasticity \cite{RewXiao}, the approach based on the multiplicative split
is gaining considerable popularity in the modeling community. This is due to numerous advantages
of the approach: Utilizing the multiplicative split, one may build a model which is
thermodynamically consistent, objective, does not exhibit any non-physical energy dissipation in the elastic range and
shear oscillation under simple shear. Moreover, the multiplicative split
is considered to be the most suitable framework for the elastoplastic analysis of single crystals \cite{Lubarda}.\footnote{
Aside from materials with a crystalline structure, the multiplicative split advocated in
this paper yields good results even for rubber-like materials \cite{Nedjar} and
geomaterials \cite{SimoMeschke}.}
Furthermore, the classical small strain models
of metal plasticity exhibit the so-called weak invariance of the solution \cite{ShutovWeakInvar}. The use of the multiplicative
split allows one to retain this weak invariance property during generalization to finite strains \cite{ShutovWeakInvar}.

Apart from the modeling of nonlinear kinematic hardening,
the extensions of the Shutov and Krei\ss ig model are used in analysis of nonlinear distortional hardening \cite{ShutovPaKr},
thermoplasticity \cite{ShutovThermal}, ductile damage \cite{ShutovDuctDamage}, and ratchetting \cite{ZHU2013}.
Moreover, the model of Shutov and Krei\ss ig was coupled with a microscopic model to account for the evolution
of metal microstructure in \cite{SilberDisloc}.
Due to the manifold applications of this model,
its numerical treatment becomes an important task. The main goal of the current paper is
to report a new efficient and robust stress integration
algorithm, which benefits from the special structure of the governing equations.

A similar phenomenological model with
a nested multiplicative split was presented in \cite{Vladimirov2008}
by Vladimirov et al.
This model was applied to many practical problems (see, among others, \cite{Vladimirov2009, Brepols2014,FrischkornReese}).
Moreover, the idea of the nested multiplicative split was also adopted by other scholars
to nonlinear kinematic hardening in \cite{DettRes, TsakmakisWulleweit, Hartmann2008, HennanAnand} and to
shape memory alloys in \cite{Helm1, Helm3, ReeseChrist, FrischkornReese}.
Due to the similarity of the modeling assumptions, the algorithm proposed in the current study
can be applied to these constitutive equations as well.

For the implementation in a displacement-based FEM, a numerical algorithm is needed which computes the stresses
at each point of Gau\ss \ integration as a function of the local deformation history \cite{SimoHughes}.
For the model of Shutov and Krei\ss ig (and its further extensions), the stress computation
is based on the integration of the underlying constitutive equations.
As usual for elasto-plasticity of metals,
the corresponding system of constitutive equations is stiff. Thus, in general, an explicit
time integration yields unacceptable results,\footnote{For the discussion of various explicit and semi-explicit
integrators, the reader is referred to \cite{Safaei} and references cited therein.}
and implicit time stepping is required to obtain a stable integration procedure.
In the current study, the preference is given to the Euler Backward
discretization of the evolution equations \cite{SimoHughes, SimMieh}.
Alternatively, one may consider its modification based on the
exponential mapping \cite{WebAnan, MiStei, Simo}.
A special class of time-stepping algorithms
is obtained by the minimization of the integrated stress power. Such a
variational approach is rather general so it can be
applied to various models with nonlinear kinematic hardening \cite{Mosler2010}.

As a result of the implicit time discretization, a system of nonlinear algebraic equations is obtained.
Unfortunately, in general, the closed-form solution for this system is not available. Therefore, a time-consuming local iterative procedure
has to be implemented to solve this system of equations \cite{ShutovKrKoo, Vladimirov2008, ReeseChrist, Brepols2014, FrischkornReese}.
In some cases, the overall system of algebraic equations can be reduced to the solution of a
scalar equation with respect to a single unknown. Such ``one-equation integrators"
are already available for many models of elastoplasticity with nonlinear kinematic hardening:
An efficient stress algorithm was constructed for a model which is based on
the assumption of small elastic strains in \cite{Luhrs}.
In \cite{Hartmann}, another anisotropic plasticity model is analyzed, which is based on
the additive decomposition of the linearized Green strain.
The authors benefit from the simple structure of the model, which resembles the infinitesimal
elasto-plasticity. 
In \cite{Zhu2014}, a problem-adapted algorithm was proposed for a
model of hypoelasto-plasticity with the additive decomposition of the strain rate.
Next, a structure-exploiting scheme for a model
with nonlinear kinematic hardening was recently proposed in
\cite{Rothe2015} in the small-strain context. In \cite{Montans2012}, an anisotropic model is considered
employing logarithmic (Hencky) strain measures and quadratic stored energy functions.
These assumptions simplify the algorithmic treatment,
since an implicit integration by exponential mapping can be carried out using small strain procedures.

At the same time, less progress has been achieved in the area of multiplicative elasto-plasticity
or multiplicative viscoelasticity coupled to hyperelasticity. Due to the nonlinear kinematics, an efficient
time-stepping becomes non-trivial even in the isotropic case.
The construction of problem-adapted integrators is
still a topic of ongoing research: An explicit update formula was
proposed in \cite{ShutovLandgraf2013} for the special case of neo-Hookean free energy.
Its generalization to elasto-plastic material with Yeoh elasticity is discussed in \cite{LandgrShuPAMM};
an application to a temperature-dependent constitutive behavior is presented in \cite{Johlitz}.
A semi-implicit scheme was proposed in \cite{HennanAnand} for a
model with nonlinear kinematic hardening implementing the Lion approach.
This scheme is constructed assuming that the elastic strains are always
coaxial with the plastic strain increment,\footnote{
This coaxiality is discussed after equation (12.11) in \cite{HennanAnand}.} which is a
very restrictive assumption in case of plastic anisotropy.
Moreover, in order to simplify the algorithmic treatment, the authors assume in \cite{HennanAnand} that the plastic strain increments are small.
Another hybrid explicit/implicit procedure was proposed in \cite{SilbermannPAMM}, which is more robust than a purely explicit
procedure.

In this paper,
the overall system of equations is subdivided (partitioned) into
coupled subsystems.
For each
of the subsystems, \emph{simple closed-form solutions are obtained}. An efficient
time stepping procedure for the fully coupled system of
algebraic equations is then constructed by a tailored combination of these explicit solutions.
As a result, a new algorithm is obtained, which is called partitioned Euler Backward method (PEBM).
Within PEBM, only a scalar equation with respect to a single unknown
inelastic strain increment has to be solved iteratively.

In the current paper, the accuracy of the stress computation by the proposed PEBM is tested numerically using
material parameters of 5754-O aluminum alloy and 42CrMo4 steel. Concerning the accuracy, the new algorithm is
comparable to the Euler Backward method with subsequent correction of incompressibility (EBMSC) and
the classical exponential method (EM). The novel PEBM is superior to the conventional EBMSC and EM concerning the
computational efficiency:
For a model with $N$ backstress tensors,
the EBMSC and EM require a solution of $6 N+7$ nonlinear scalar equations at each Gau\ss \ point.
At the same time, the PEBM requires a solution of a scalar equation with respect
to a single unknown inelastic strain increment. In particular, the introduction of additional
backstresses is possible almost without any increase in computational costs.

We close this introduction with a few words regarding notation.
Throughout this paper, bold-faced symbols denote first- and second-rank tensors in $\mathbb{R}^3$.
A coordinate-free tensor formalism is used in this work;
$\mathbf{1}$ stands for the second-rank identity tensor;
the deviatoric part of a tensor is denoted by
$\mathbf A^{\text{D}} := \mathbf A - \frac{1}{3} \text{tr}(\mathbf A) \mathbf 1$, where
$\text{tr}(\mathbf A)$ stands for the trace.
The overline $\overline{(\cdot)}$ stands for the unimodular
part of a tensor
\begin{equation}\label{BarDef}
\overline{\mathbf{A}}=(\det \mathbf{A})^{-1/3} \mathbf{A}.
\end{equation}
In description of algorithms, the notation $\mathbf{A} \gets \mathbf{B}$ means a computation step
where the variable $\mathbf{A}$ obtains the new value $\mathbf{B}$.

\section{MATERIAL MODEL}

\subsection{System of constitutive equations}

First, we recall the viscoplastic material model proposed by Shutov and Krei\ss ig in \cite{ShutovKrVisc}.
Here we consider its modification with two Armstrong-Frederick backstresses (cf. \cite{ShutovKuprin}). A generalization of the algorithm to
cover an arbitrary number of backstresses is straightforward.
For simplicity of the numerical implementation, we adopt the Lagrangian (material) formulation
of the constitutive equations.\footnote{An alternative derivation of the model is presented in \cite{Broec}.}
The description of finite strain kinematics is based on the original work of Lion \cite{LionIJP}. Along with the well-known right Cauchy-Green tensor
${\mathbf C}$, we introduce
tensor-valued internal variables:
the inelastic right Cauchy-Green tensor
 ${\mathbf C}_{\text{i}}$ and two
inelastic right Cauchy-Green tensors of substructure:
 ${\mathbf C}_{\text{1i}}$ and ${\mathbf C}_{\text{2i}}$.
In order to capture the isotropic hardening we will need
the inelastic arc length (Odqvist parameter)
$s$ and its dissipative part $s_d$.

For the free energy per unit mass we assume the additive decomposition as follows
\begin{equation}\label{EnergStorage}
\psi = \psi_{\text{el}} (\mathbf C {\mathbf C_{\text{i}}}^{-1}) +
\psi_{\text{kin 1}}(\mathbf C_{\text{i}} {\mathbf C_{\text{1i}}}^{-1}) +
\psi_{\text{kin 2}}(\mathbf C_{\text{i}} {\mathbf C_{\text{2i}}}^{-1}) + \psi_{\text{iso}}(s-s_d).
\end{equation}
Here, $\psi_{\text{el}}$ stands for the energy storage due to macroscopic elastic deformations,
$\psi_{\text{kin 1}}$ and $\psi_{\text{kin 2}}$ represent two different storage mechanisms associated to the
kinematic hardening, and $\psi_{\text{iso}}$ is the energy storage
related to the isotropic hardening.\footnote{Some additional
energy-storage terms, which are not
related to any hardening mechanism, can be added to the right-hand side of \eqref{EnergStorage}.
Such terms can play an important role in thermoplasticity \cite{ShutovThermal}.}
The functions $\psi_{\text{el}}$, $\psi_{\text{kin 1}}$ and $\psi_{\text{kin 2}}$ are assumed to be isotropic.
Let $\tilde{\mathbf T}$ and $\tilde{\mathbf X}$ denote respectively the 2nd Piola-Kirchhoff stress tensor and
the total backstress tensor operating on the reference configuration.
The backstress describes the translation of the yield surface in the stress space, and it
will be decomposed into two partial
backstresses $\tilde{\mathbf X}_{1}$ and $\tilde{\mathbf X}_{2}$.
The isotropic expansion of the yield surface is captured by the scalar quantity $R$.
The rate-dependent overstress $f$ is introduced to account for the viscous effects.

In order to formulate the system of constitutive
equations, consider a local initial value problem.
The local deformation history ${\mathbf C}(t)$ is assumed to be
known in the time interval $t \in [0,T]$.
The corresponding stress response is governed by
the following system of ordinary differential and algebraic equations
\begin{equation}\label{prob1}
\dot{\mathbf C}_{\text{i}} = 2 \frac{\displaystyle
\lambda_{\text{i}}}{\displaystyle \mathfrak{F}}
\big( \mathbf C \tilde{\mathbf T} -
\mathbf C_{\text{i}} \tilde{\mathbf X} \big)^{\text{D}} \mathbf C_{\text{i}}, \quad
\mathbf C_{\text{i}}|_{t=0} = \mathbf C_{\text{i}}^0, \
\det \mathbf C_{\text{i}}^0 =1, \ \mathbf C_{\text{i}}^0 = (\mathbf C_{\text{i}}^0)^{\text{T}} > 0,
\end{equation}
\vspace{-5mm}
\begin{equation}\label{prob21}
\dot{\mathbf C}_{\text{1i}} =
2 \lambda_{\text{i}} \varkappa_1 (\mathbf C_{\text{i}}
\tilde{\mathbf X_{1}} \big)^{\text{D}} \mathbf C_{\text{1i}}, \quad
\mathbf C_{\text{1i}}|_{t=0} = \mathbf C_{\text{1i}}^0, \
\det \mathbf C_{\text{1i}}^0 =1, \ \mathbf C_{\text{1i}}^0 = (\mathbf C_{\text{1i}}^0)^{\text{T}} > 0,
\end{equation}
\begin{equation}\label{prob22}
\dot{\mathbf C}_{\text{2i}} =
2 \lambda_{\text{i}} \varkappa_2 (\mathbf C_{\text{i}}
\tilde{\mathbf X_{2}} \big)^{\text{D}} \mathbf C_{\text{2i}}, \quad
\mathbf C_{\text{2i}}|_{t=0} = \mathbf C_{\text{2i}}^0, \
\det \mathbf C_{\text{2i}}^0 =1, \ \mathbf C_{\text{2i}}^0 = (\mathbf C_{\text{2i}}^0)^{\text{T}} > 0,
\end{equation}
\begin{equation}\label{prob3}
\dot{s}= \sqrt{\frac{\displaystyle 2}{\displaystyle 3}} \lambda_{\text{i}}, \quad
\dot{s}_{\text{d}}= \frac{\displaystyle \beta}{\displaystyle \gamma} \dot{s} R,  \quad
s|_{t=0} = s^0, \ s_{\text{d}}|_{t=0} = s_{\text{d}}^0,
\end{equation}
\begin{equation}\label{prob4}
\tilde{\mathbf T} =
2 \rho_{\scriptscriptstyle \text{R}}
\frac{\displaystyle \partial \psi_{\text{el}}
(\mathbf C {\mathbf C_{\text{i}}}^{-1})}
{\displaystyle \partial
\mathbf{C}}\big|_{\mathbf C_{\text{i}} = \text{const}},
\quad \tilde{\mathbf X} = \tilde{\mathbf X}_{1} + \tilde{\mathbf X}_{2},
\end{equation}
\begin{equation}\label{prob42}
\tilde{\mathbf X}_{1} =
2 \rho_{\scriptscriptstyle \text{R}}
\frac{\displaystyle \partial \psi_{\text{kin 1}}(\mathbf C_{\text{i}} {\mathbf C_{\text{1i}}}^{-1})}
{\displaystyle \partial \mathbf C_{\text{i}}}\big|_{\mathbf C_{\text{1i}} = \text{const}}, \quad
\tilde{\mathbf X}_{2} =
2 \rho_{\scriptscriptstyle \text{R}}
\frac{\displaystyle \partial \psi_{\text{kin 2}}(\mathbf C_{\text{i}} {\mathbf C_{\text{2i}}}^{-1})}
{\displaystyle \partial \mathbf C_{\text{i}}}\big|_{\mathbf C_{\text{2i}} = \text{const}},
\end{equation}
\begin{equation}\label{prob5}
 R = \rho_{\scriptscriptstyle \text{R}}
 \frac{\displaystyle \partial \psi_{\text{iso}}(s-s_d)}{\displaystyle \partial s},
\end{equation}
\begin{equation}\label{prob6}
\lambda_{\text{i}}= \frac{\displaystyle 1}{\displaystyle
\eta}\Big\langle \frac{\displaystyle 1}{\displaystyle f_0}
f \Big\rangle^{m}, \quad
f= \mathfrak{F}- \sqrt{\frac{2}{3}} \big[ K + R \big], \quad
\mathfrak{F}= \sqrt{\text{tr}
\big[ \big( \mathbf C \tilde{\mathbf T} -
\mathbf C_{\text{i}} \tilde{\mathbf X} \big)^{\text{D}} \big]^2 }.
\end{equation}
Here,
$\rho_{\scriptscriptstyle \text{R}} > 0$,
$\varkappa_1 \geq 0$, $\varkappa_2 \geq 0$, $\beta \geq 0$, $\gamma \in \mathbb{R}$, $\eta \geq 0$, $m \geq 1$, $K > 0$
are constant material parameters,
the real-valued constitutive functions $\psi_{\text{el}}$, $\psi_{\text{kin 1}}$, $\psi_{\text{kin 2}}$, $\psi_{\text{iso}}$ are assumed to be known.
The constant $f_0 = 1$ MPa is adopted to obtain a non-dimensional term in the bracket; thus, $f_0$ is not a material parameter.
The symbol $\lambda_{\text{i}}(t)$ stands for the inelastic multiplier; $\mathfrak{F}(t)$ denotes the norm of the driving force;
the superposed dot denotes the material time derivative: $\frac{d}{d t} \mathbf A = \dot{\mathbf A}$;
$(\cdot)^{\text{D}}$ is the deviatoric part of a second-rank tensor;
$\langle x \rangle := \max (x,0)$ is the Macaulay bracket.

For $\psi_{\text{iso}}$ we assume
\begin{equation}\label{isotrAssum}
\rho_{\scriptscriptstyle \text{R}}  \psi_{\text{iso}}(s-s_{\text{d}}) = \frac{\gamma}{2} (s - s_{\text{d}})^2,
\end{equation}
which yields
\begin{equation}\label{isotrAssum2}
R = \gamma (s - s_{\text{d}}).
\end{equation}
Note that this restrictive assumption for $\psi_{\text{iso}}$ is not essential for the
efficient time stepping scheme, which will be presented in the following.
For the remaining part of the energy storage we adopt potentials of compressible neo-Hookean material
\begin{equation}\label{spec11}
\rho_{\scriptscriptstyle \text{R}}  \psi_{\text{el}}(\mathbf{C} \mathbf{C}_{\text{i}}^{-1})=
\frac{k}{50}\big( (\text{det} \mathbf{C} \mathbf{C}_{\text{i}}^{-1})^{5/2} + (\text{det} \mathbf{C} \mathbf{C}_{\text{i}}^{-1})^{-5/2} -2 \big)+
\frac{\mu}{2} \big( \text{tr} \overline{\mathbf{C} \mathbf{C}_{\text{i}}^{-1}} - 3 \big),
\end{equation}
\begin{equation}\label{spec12}
\rho_{\scriptscriptstyle \text{R}}
\psi_{\text{kin 1}}(\mathbf{C}_{\text{i}} \mathbf{C}_{\text{i1}}^{-1})=
\frac{c_1}{4}\big( \text{tr} \overline{\mathbf{C}_{\text{i}} \mathbf{C}_{\text{1i}}^{-1}} - 3 \big), \quad
\rho_{\scriptscriptstyle \text{R}}
\psi_{\text{kin 2}}(\mathbf{C}_{\text{i}} \mathbf{C}_{\text{2i}}^{-1})=
\frac{c_1}{4}\big( \text{tr} \overline{\mathbf{C}_{\text{i}} \mathbf{C}_{\text{2i}}^{-1}} - 3 \big),
\end{equation}
where $k >0$, $\mu > 0$, $c1 > 0$, and $c2 > 0$ are material parameters.

\textbf{Remark 1.} The volumetric part of the strain energy function \eqref{spec11} was
proposed by Hartmann and Neff in \cite{HartmannNeff}.
It has an important advantage over many alternatives since it provides a convex function of $\det \mathbf{F}$. Moreover,
the energy and the hydrostatic stress component tend to infinity as $\det \mathbf{F} \rightarrow 0$ or $\det \mathbf{F} \rightarrow \infty$. Thus, ansatz \eqref{spec11} predicts a
physically reasonable volumetric stress response even for finite elastic volume changes.
Note that the algorithm, presented in the current study, can be applied irrespective of the specific choice of
the volumetric strain energy function. $ \Box $

Using \eqref{spec11}, \eqref{spec12} and the incompressibility condition $\text{det} \mathbf C_{\text{i}} =1$ which
will be discussed in the following, relations \eqref{prob4} and \eqref{prob42} are rewritten as
\begin{equation}\label{trans41}
\tilde{\mathbf T} = \frac{\displaystyle k}{\displaystyle 10} \
\big( (\text{det} \mathbf C)^{5/2}-(\text{det} \mathbf C)^{-5/2} \big) \
\mathbf C^{-1} + \mu \ \mathbf C^{-1} (\overline{\mathbf C}
\mathbf C_{\text{i}}^{-1})^{\text{D}},
\end{equation}
\begin{equation}\label{trans42}
\tilde{\mathbf X}_1 = \frac{c_1}{2} \ \mathbf C_{\text{i}}^{-1}
(\mathbf C_{\text{i}} \mathbf C_{\text{1i}}^{-1})^{\text{D}}, \quad
\tilde{\mathbf X}_2 = \frac{c_2}{2} \ \mathbf C_{\text{i}}^{-1}
(\mathbf C_{\text{i}} \mathbf C_{\text{2i}}^{-1})^{\text{D}}.
\end{equation}
Substituting these results into $\eqref{prob1}$ --- $\eqref{prob22}$,
we obtain the reduced evolution equations
\begin{equation}\label{ReducedEvolut}
\dot{\mathbf C}_{\text{i}} = 2 \frac{\displaystyle
\lambda_{\text{i}}}{\displaystyle \mathfrak{F}}
 \big( \mu \ (\overline{\mathbf C}
\mathbf C_{\text{i}}^{-1})^{\text{D}}  -
 \frac{c_1}{2} \ (\mathbf C_{\text{i}} \mathbf C_{\text{1i}}^{-1})^{\text{D}} -
 \frac{c_2}{2} \ (\mathbf C_{\text{i}} \mathbf C_{\text{2i}}^{-1})^{\text{D}}  \big)  \mathbf C_{\text{i}},
 \end{equation}
\begin{equation}\label{ReducedEvolut2}
 \dot{\mathbf C}_{\text{1i}} =
\lambda_{\text{i}} \ \varkappa_1 \ c_1 \ (\mathbf C_{\text{i}} \mathbf C_{\text{1i}}^{-1})^{\text{D}} \mathbf C_{\text{1i}}, \quad
\dot{\mathbf C}_{\text{2i}} =
\lambda_{\text{i}} \ \varkappa_2 \ c_2 \ (\mathbf C_{\text{i}} \mathbf C_{\text{2i}}^{-1})^{\text{D}} \mathbf C_{\text{2i}}.
\end{equation}

\subsection{Properties of the material model}

The material model is thermodynamically consistent.\footnote{
The thermodynamic consistency of the model was proved in \cite{ShutovKrVisc} for the special case of only one backstress.
This proof can be easily generalized to cover numerous backstress tensors.} It combines the nonlinear hardening of Armstrong-Frederick type and the
isotropic hardening of Voce type. The exact solution of \eqref{prob1} -- \eqref{prob6} exhibits the following geometric property
\begin{equation}\label{geopr0}
\mathbf{C}_{\text{i}}, \mathbf{C}_{\text{1i}}, \mathbf{C}_{\text{2i}} \in M, \quad \text{where} \
M := \big\{ \mathbf B \in Sym: \text{det} \mathbf B =1 \big\}.
\end{equation}
Thus, \eqref{prob1} -- \eqref{prob6} is a system of differential and
algebraic equations on the manifold. The restriction $\text{det} (\mathbf{C}_{\text{i}}) =1$
corresponds to the incompressibility of the inelastic flow, typically assumed for metals;
the conditions $\text{det} (\mathbf{C}_{\text{1i}}) =1$ and  $\text{det} (\mathbf{C}_{\text{2i}}) =1$
reflect the incompressibility on the substructural level, which is a modeling assumption.
As was shown in \cite{ShutovKrStab},
these incompressibility conditions should be exactly satisfied by the numerical algorithm in order to suppress
the accumulation of the integration error.
Another important restriction is as follows: For the exact solution, the
tensor-valued variables $\mathbf{C}_{\text{i}}$, $\mathbf{C}_{\text{1i}}$ and $\mathbf{C}_{\text{2i}}$
are positive definite. Also from the physical standpoint, these quantities must remain positive definite, since they
represent some inelastic metric tensors of Cauchy-Green type.

One important property of the material model, which will play
a crucial role in the current study, is the \emph{weak invariance}
under arbitrary isochoric changes of the reference configuration.
Generally speaking, the weak invariance corresponds to a generalized (weak) symmetry of the material.
Like any other symmetry property, the weak invariance should be inhereted by the discretized constitutive equations.
For the general definition of the weak invariance the reader is referred to \cite{ShutovWeakInvar}.
For the concrete material model analyzed in the current study, the weak invariance is formulated as follows.
Let $\textbf{F}_0$ be arbitrary second rank tensor such that $\text{det}(\textbf{F}_0) = 1$.
If the prescribed local loading $\textbf{C}(t), \ t \in [0,T]$ is replaced by a new loading history
$\textbf{C}^{\text{new}}(t) := \textbf{F}^{-\text{T}}_0 \  \textbf{C}(t) \ \textbf{F}^{-1}_0$ and the
initial conditions are transformed according to
\begin{equation*}\label{IniCond}
\textbf{C}_{\text{i}}^{\text{new}}|_{t=0} = \textbf{F}^{-\text{T}}_0 \ \textbf{C}_{\text{i}}|_{t=0} \ \textbf{F}^{-1}_0, \quad
\textbf{C}_{\text{1i}}^{\text{new}}|_{t=0} = \textbf{F}^{-\text{T}}_0 \ \textbf{C}_{\text{1i}}|_{t=0} \ \textbf{F}^{-1}_0, \quad
\textbf{C}_{\text{2i}}^{\text{new}}|_{t=0} = \textbf{F}^{-\text{T}}_0 \ \textbf{C}_{\text{2i}}|_{t=0} \ \textbf{F}^{-1}_0,
\end{equation*}
then the second Piola-Kirchhoff tensor, predicted by the model, is transformed according to
\begin{equation}\label{WeakInvariance_Reference}
\widetilde{\mathbf{T}}^{\text{new}} (t)  =
\mathbf{F}_0 \ \widetilde{\mathbf{T}} (t) \
\mathbf{F}_0^{\text{T}} \quad \text{for} \ t \in [0,T].
\end{equation}
Moreover, if the history of the deformation gradient is replaced through
$\textbf{F}^{\text{new}}(t) :=  \textbf{F}(t) \ \textbf{F}^{-1}_0$,
then the same Cauchy stresses $\textbf{T}$ will be predicted by the model
\begin{equation}\label{WeakInvariance}
\textbf{T}^{\text{new}}(t) = \textbf{T}(t) \quad \text{for} \ t \in [0,T].
\end{equation}
For this invariance property, it is necessary and sufficient that the internal variables are transformed according
to the following transformation rules (cf. \cite{ShutovPfeIhl})
\begin{equation}\label{WeakInvariance_InterVariables}
\textbf{C}_{\text{i}}^{\text{new}}(t) = \textbf{F}^{-\text{T}}_0 \ \textbf{C}_{\text{i}}(t) \ \textbf{F}^{-1}_0, \quad
\textbf{C}_{\text{1i}}^{\text{new}}(t) = \textbf{F}^{-\text{T}}_0 \ \textbf{C}_{\text{1i}}(t) \ \textbf{F}^{-1}_0, \quad
\textbf{C}_{\text{2i}}^{\text{new}}(t) = \textbf{F}^{-\text{T}}_0 \ \textbf{C}_{\text{2i}}(t) \ \textbf{F}^{-1}_0.
\end{equation}

In metals, the plastic flow may exhibit a fluid-like feature that there is no any preferred configuration.
This observation motivates the notion of the weak invariance from the physical standpoint.
For a more detailed discussion, the
reader is referred to \cite{ShutovWeakInvar}.
Some of the commonly used approaches to finite strain plasticity are compatible with the weak invariance requirement, but some of them not \cite{ShutovWeakInvar}.
The weak invariance brings the advantage that the constitutive equations are independent of
the choice of the reference configuration. In particular, this property simplifies the
analysis of multi-stage forming operations \cite{ShutovPfeIhl}.
In the current study, the weak invariance is utilized as a guide in constructing numerical algorithms.

\section{Implicit time stepping}

Let us consider a time step $t_n \mapsto t_{n+1}$ with $\Delta t := t_{n+1} - t_n >0$.
Within a local setting, we may suppose that ${}^{n+1} \mathbf C := \mathbf C (t_{n+1})$ and ${}^{n} \mathbf C := \mathbf C (t_{n})$
are known and the internal variables are given at $t_n$ by
${}^{n} \mathbf C_{\text{i}}$, ${}^{n} \mathbf C_{\text{1i}}$, ${}^{n} \mathbf C_{\text{2i}}$, ${}^{n} s$, ${}^{n} s_{\text{d}}$.
In order to compute the stress tensor
${}^{n+1} \tilde{\mathbf T} $ we need to update the internal variables.

For what follows, it is useful to introduce the incremental inelastic strain
\begin{equation}\label{ineinc}
\xi := \Delta t \ {}^{n+1} \lambda_{\text{i}} \geq 0,
\end{equation}
which is a non-dimensional quantity. For most practical applications,
we may assume that $ \xi \leq 0.2$.\footnote{This upper bound is based on experience (see also \cite{Hartmann}).
Since the material response is highly path dependent, one should avoid
plastic increments larger than 20\%.}
A new implicit numerical method will be presented in this section which is
based on some closed-form solutions for elementary decoupled problems.
For this new method, the time stepping procedure will be reduced
to the solution of a scalar equation
with respect to the unknown $\xi$.

\subsection{Update of $s$ and $s_\text{d}$ for a given $\xi$}

Integrating $\eqref{prob3}_1$ and $\eqref{prob3}_2$ from $t_n$ to $t_{n+1}$ and taking
\eqref{isotrAssum2} into account, we obtain
the isotropic hardening
${}^{n+1} R$, the arc-length ${}^{n+1} s$, and its dissipative part ${}^{n+1} s_{\text{d}}$ as explicit functions of $\xi$
\begin{equation}\label{updateR}
{}^{n+1} R (\xi) = \frac{{}^{\text t} R + \sqrt{2/3} \gamma \xi}{1 + \sqrt{2/3} \beta \xi},
\quad \text{where} \quad  {}^{\text t} R := \gamma({}^{n} s - {}^{n} s_{\text{d}}),
\end{equation}
\begin{equation}\label{updateSSd}
{}^{n+1} s (\xi) = {}^{n} s + \sqrt{2/3} \ \xi, \quad
{}^{n+1} s_{\text{d}} (\xi) = {}^{n} s_{\text{d}} + \frac{\beta}{\gamma} \ \sqrt{2/3} \ \xi \ {}^{n+1} R (\xi).
\end{equation}

\textbf{Remark 2.} In case of constant $\gamma$ and $\beta$, a closed-form solution can be derived for the system
\eqref{prob3}, \eqref{isotrAssum2}, which can also be used for the time integration.  $\Box$

\subsection{Update of $\textbf{C}_{\text{i}}$ for the known  ${}^{n+1}\textbf{C}_{\text{1i}}$, ${}^{n+1}\textbf{C}_{\text{2i}}$, $\xi$}

Assume for a moment that ${}^{n+1}\textbf{C}_{\text{1i}}$, ${}^{n+1}\textbf{C}_{\text{2i}}$, and $\xi > 0$ are known.
In this section we update $\textbf{C}_{\text{i}}$.
As a first step, we express the norm of the driving force as a function of $\xi$:
Combining relations $\eqref{prob6}_1$ and $\eqref{prob6}_2$ we arrive at
\begin{equation}\label{F2}
\mathfrak{F} = \mathfrak{F}_2(\xi) := f_0 \big(\frac{\eta \xi}{\Delta t}\big)^{1/m} + \sqrt{2/3}(K + \ {}^{n+1}R(\xi)).
\end{equation}
Next, we rewrite the evolution equation \eqref{ReducedEvolut} as follows
\begin{equation}\label{ReducedEvolutCi2}
\dot{\mathbf C}_{\text{i}} = 2 \frac{\displaystyle
\lambda_{\text{i}}}{\displaystyle \mathfrak{F}}
 \big( \mu \ \overline{\mathbf C}
 -  \frac{c_1}{2} \ \mathbf C_{\text{i}} \mathbf C_{\text{1i}}^{-1} \mathbf C_{\text{i}}
 -  \frac{c_2}{2} \ \mathbf C_{\text{i}} \mathbf C_{\text{2i}}^{-1} \mathbf C_{\text{i}} \big) + \beta \mathbf C_{\text{i}},
\end{equation}
where
\begin{equation}\label{EqForBeta}
\beta = -\frac{2}{3} \frac{\displaystyle
\lambda_{\text{i}}}{\displaystyle \mathfrak{F}}\big( \mu \ \text{tr}(\overline{\mathbf C} \mathbf C_{\text{i}}^{-1})
- \frac{c_1}{2} \ \text{tr}(\mathbf C_{\text{i}} \mathbf C_{\text{1i}}^{-1})
- \frac{c_2}{2} \ \text{tr}(\mathbf C_{\text{i}} \mathbf C_{\text{2i}}^{-1})\big) \in \mathbb{R}.
\end{equation}
Implementing the Euler Backward method (EBM) to \eqref{ReducedEvolutCi2}, and replacing
$\mathfrak{F}$ through  $\mathfrak{F}_2(\xi)$ we obtain
\begin{equation}\label{ReducedEvolut22}
{}^{n+1} \mathbf C_{\text{i}} = {}^{n} \mathbf C_{\text{i}} +  2 \frac{\displaystyle
\xi }{\displaystyle \mathfrak{F}_2}
 \big( \mu \ {}^{n+1} \overline{\mathbf C}
 -  \frac{c_1}{2} \ {}^{n+1} \mathbf C_{\text{i}} \ {}^{n+1}\mathbf C_{\text{1i}}^{-1} \
 {}^{n+1} \mathbf C_{\text{i}}
 -  \frac{c_2}{2} \ {}^{n+1} \mathbf C_{\text{i}} \ {}^{n+1}\mathbf C_{\text{2i}}^{-1} \
 {}^{n+1} \mathbf C_{\text{i}}  \big) + \beta \Delta t \ {}^{n+1} \mathbf C_{\text{i}},
\end{equation}
where the scalar $\beta$ is unknown since it depends on the unknown ${}^{n+1}\textbf{C}_{\text{i}}$.

Unfortunately, due to the additive nature of the classical EBM, this method
violates the incompressibility restriction $\text{det} (\mathbf C_{\text{i}}) =1$.
In order to enforce the exact incompressibility, we may modify the right-hand side of \eqref{ReducedEvolut22} by
introducing additional (correction) term $\varepsilon \mathbf{P}$, thus arriving at
\begin{equation}\label{ReducedEvolut22mod}
{}^{n+1} \mathbf C_{\text{i}} = {}^{n} \mathbf C_{\text{i}} +  2 \frac{\displaystyle
\xi }{\displaystyle \mathfrak{F}_2}
 \big( \mu \ {}^{n+1} \overline{\mathbf C}
 -  \frac{c_1}{2} \ {}^{n+1} \mathbf C_{\text{i}} \ {}^{n+1}\mathbf C_{\text{1i}}^{-1} \ {}^{n+1} \mathbf C_{\text{i}}
 -  \frac{c_2}{2} \ {}^{n+1} \mathbf C_{\text{i}} \ {}^{n+1}\mathbf C_{\text{2i}}^{-1} \
 {}^{n+1} \mathbf C_{\text{i}}  \big) + \beta \Delta t \ {}^{n+1} \mathbf C_{\text{i}} + \varepsilon \mathbf{P}.
\end{equation}
Here, $\varepsilon$ is a suitable scalar and $\mathbf{P}$ is a suitable second-rank tensor.\footnote{Other modification of EBM used
to enforce the incompressibility were presented in \cite{Helm2, Vladimirov2008}.}
The scalar $\varepsilon$
is determined such that $\text{det} ({}^{n+1} \mathbf C_{\text{i}}) =1$.
In some publications, the authors suggest the correction
of diagonal terms in ${}^{n+1} \mathbf C_{\text{i}}$ to enforce the incompressibility. Such a correction
corresponds to $\mathbf{P} =\mathbf{1}$, but this choice of $\mathbf{P}$ violates
the weak invariance of the solution (see Appendix A). In the current study we put $\mathbf{P} = {}^{n+1} \mathbf C_{\text{i}}$.
Such a choice is \emph{compatible with the weak invariance} (see Appendix A).
After some algebraic manipulations we obtain from \eqref{ReducedEvolut22mod}
\begin{equation}\label{ReducedEvolut24}
z \ {}^{n+1} \mathbf C_{\text{i}} = {}^{n} \mathbf C_{\text{i}} +  2 \frac{\displaystyle
\xi }{\displaystyle \mathfrak{F}_2}
 \big( \mu \ {}^{n+1} \overline{\mathbf C}
 -  \frac{c_1}{2} \ {}^{n+1} \mathbf C_{\text{i}} \ {}^{n+1}\mathbf C_{\text{1i}}^{-1}
 {}^{n+1} \mathbf C_{\text{i}}
 -  \frac{c_2}{2} \ {}^{n+1} \mathbf C_{\text{i}} \ {}^{n+1}\mathbf C_{\text{2i}}^{-1} \
 {}^{n+1} \mathbf C_{\text{i}}  \big),
\end{equation}
where the scalar unknown $z := 1 - \beta \Delta t - \varepsilon$ is determined from the
incompressibility condition $\text{det} ({}^{n+1} \mathbf C_{\text{i}}) =1$.
Next, we introduce the following abbreviations
\begin{equation}\label{AbbrevPhiC}
c:= c_1 + c_2, \quad \mathbf{\Phi} := \frac{\displaystyle c_1}{c} \ {}^{n+1}\mathbf C_{\text{1i}}^{-1} +
\frac{\displaystyle c_2}{c} \ {}^{n+1}\mathbf C_{\text{2i}}^{-1}.
\end{equation}
Using these, we rewrite \eqref{ReducedEvolut24} as follows
\begin{equation}\label{ReducedEvolut3}
z \ {}^{n+1} \mathbf C_{\text{i}} = {}^{n} \mathbf C_{\text{i}} +  2 \frac{\displaystyle
\xi }{\displaystyle \mathfrak{F}_2}
\mu \ {}^{n+1} \overline{\mathbf C}
-  \frac{\displaystyle \xi c}{\displaystyle \mathfrak{F}_2} \
{}^{n+1} \mathbf C_{\text{i}} \ \mathbf{\Phi} \  {}^{n+1} \mathbf C_{\text{i}}.
\end{equation}
Further, multiplying both sides of \eqref{ReducedEvolut3} with $\mathbf{\Phi}^{1/2}$ from the left and right, we obtain
\begin{equation}\label{ReducedEvolut4}
z \ \mathbf{\Phi}^{1/2} \ {}^{n+1} \mathbf C_{\text{i}} \ \mathbf{\Phi}^{1/2} = \mathbf{\Phi}^{1/2} \big[ {}^{n} \mathbf C_{\text{i}}  + 2 \frac{\displaystyle
\xi }{\displaystyle \mathfrak{F}_2}
 \mu \ {}^{n+1} \overline{\mathbf C} \ \big] \mathbf{\Phi}^{1/2} -\frac{\displaystyle
\xi c}{\displaystyle \mathfrak{F}_2} \ \mathbf{\Phi}^{1/2}
 \ {}^{n+1} \mathbf C_{\text{i}} \ \mathbf{\Phi} \  {}^{n+1} \mathbf C_{\text{i}} \ \mathbf{\Phi}^{1/2}.
\end{equation}
Now let us introduce one more abbreviation
\begin{equation}\label{Abbreviation}
\mathbf A : = \mathbf{\Phi}^{1/2} \big[ {}^{n} \mathbf C_{\text{i}}  + 2 \frac{\displaystyle
\xi }{\displaystyle \mathfrak{F}_2}
\mu \ {}^{n+1} \overline{\mathbf C} \ \big] \mathbf{\Phi}^{1/2}.
\end{equation}
The tensor $\mathbf A$ is known. In order to simplify the system of equations, we rewrite it with respect to a new unknown variable
\begin{equation}\label{NewUnknown}
\mathbf Y := \mathbf{\Phi}^{1/2}  \ {}^{n+1} \mathbf C_{\text{i}} \ \mathbf{\Phi}^{1/2}.
\end{equation}
Substituting  relations \eqref{Abbreviation} and \eqref{NewUnknown} into \eqref{ReducedEvolut4}, we obtain a
quadratic equation with respect to $\mathbf Y$
\begin{equation}\label{QuadrEquation}
z \ \mathbf Y = \mathbf A -  \frac{\displaystyle
\xi c}{\displaystyle \mathfrak{F}_2} \mathbf Y^2.
\end{equation}
Recall that, for physical reasons, the tensor ${}^{n+1} \mathbf C_{\text{i}}$ must be positive definite.
Therefore, the physically reasonable $\mathbf Y$ is also positive definite.
Thus, the correct solution of \eqref{QuadrEquation} is given by
\begin{equation}\label{ReducedEvolut5}
\mathbf Y = \frac{\displaystyle
 \mathfrak{F}_2}{\displaystyle 2 \xi c} \Big[ -z \textbf{1} + \big(z^2 \textbf{1} +  4 \frac{\displaystyle
\xi c}{\displaystyle \mathfrak{F}_2} \textbf{A} \big)^{1/2} \Big].
\end{equation}
Unfortunately, this relation is prone to round-off errors if evaluated step-by-step (see Appendix B).
Thus, a reliable method should be used
to compute the matrix function on the right-hand side of \eqref{ReducedEvolut5}.
Using the algebraic identity $(\mathbf X^{1/2} -\mathbf 1) (\mathbf X^{1/2} + \mathbf 1) = \mathbf X - \mathbf 1$,
the closed-form solution \eqref{ReducedEvolut5} can be recast in the form
\begin{equation}\label{ReducedEvolut5notCont}
\mathbf Y = 2 \textbf{A} \Big[  \Big( z^2 \textbf{1} + \frac{4 \xi c}{\mathfrak{F}_2} \textbf{A} \Big)^{1/2} + z \textbf{1} \Big]^{-1}.
\end{equation}
In the exact arithmetics this formula is equivalent to \eqref{ReducedEvolut5} but it is more
advantageous from the numerical standpoint (see Appendix B).

Now it remains to identify the unknown parameter $z$.
Recall that $z$ can be estimated using the incompressibility condition $\det ({}^{n+1} \mathbf C_{\text{i}}) =1$,
which is equivalent to $\det (\textbf{Y}) = \det (\mathbf{\Phi})$.
As shown in Appendix C, a simple formula can be obtained using the
perturbation method for small $\frac{\displaystyle c \ \xi}{\displaystyle \mathfrak{F}_2}$:
\begin{equation}\label{EstimOfz}
z = z_0 - \frac{\text{tr} \mathbf A}{3 z_0} \frac{c \ \xi}{\mathfrak{F}_2} + O\Big(\Big(\frac{c \ \xi}{\mathfrak{F}_2}\Big)^2\Big),
\quad \text{where} \
z_0 :=  \Big(\frac{\det \mathbf A}{\det{\mathbf \Phi}}\Big)^{1/3}.
\end{equation}
Neglecting the terms $O\Big(\Big(\frac{c \ \xi}{\mathfrak{F}_2}\Big)^2\Big)$, a reliable estimation is obtained,
which is accurate for $c=0$. As will be shown in the following sections, this formula allows us to obtain an accurate
and robust local numerical procedure even for finite values of $\frac{\displaystyle c \ \xi}{\displaystyle \mathfrak{F}_2}$.

After $z$ is evaluated, $\mathbf{Y}$ is computed using \eqref{ReducedEvolut5notCont} and $\mathbf C_{\text{i}} $ is updated through
\begin{equation}\label{UpdateCi}
{}^{n+1} \mathbf C^*_{\text{i}} := \mathbf{\Phi}^{-1/2}  \ \mathbf Y \ \mathbf{\Phi}^{-1/2}.
\end{equation}
Note that this solution exactly preserves the weak invariance since
it is the solution for system \eqref{ReducedEvolut24} which itself is weakly invariant.
Next, since the variable $z$ is not computed exactly, the incompressibility condition can be violated.
Thus, to enforce the incompressibility, a final correction step is needed
\begin{equation}\label{UpdateCi_correct}
{}^{n+1} \mathbf C_{\text{i}} = \overline{{}^{n+1} \mathbf C^*_{\text{i}}}.
\end{equation}

The procedure, described in this subsection, yields ${}^{n+1} \mathbf C_{\text{i}}$ as a function of ${}^{n+1} \mathbf C_{\text{1i}}$,
${}^{n+1} \mathbf C_{\text{2i}}$, and $\xi$:
\begin{equation}\label{UpdateCi2}
{}^{n+1} \mathbf C_{\text{i}} = \mathfrak{C}_{\text{i}} ({}^{n+1} \mathbf C_{\text{1i}}, {}^{n+1} \mathbf C_{\text{2i}}, \xi).
\end{equation}

\subsection{Update of $\textbf{C}_{\text{1i}}$ and $\textbf{C}_{\text{2i}}$ for the known  ${}^{n+1} \textbf{C}_{\text{i}}$ and $\xi$}

Let us update $\textbf{C}_{\text{1i}}$ using evolution equation $\eqref{ReducedEvolut2}_1$.
Its Euler Backward discretization yields
\begin{equation}\label{UpdateCii}
{}^{n+1} \mathbf C_{\text{1i}} = {}^{n} \mathbf C_{\text{1i}} +
\xi \ \varkappa_1 \ c_1 \ ({}^{n+1} \mathbf C_{\text{i}} \  ({}^{n+1}\mathbf C_{\text{1i}})^{-1})^{\text{D}} \ {}^{n+1} \mathbf C_{\text{1i}}.
\end{equation}
Next, we rewrite this algebraic equation as follows
\begin{equation}\label{UpdateCii2}
{}^{n+1} \mathbf C_{\text{1i}} = {}^{n} \mathbf C_{\text{1i}} +
\xi \ \varkappa_1 \ c_1 \ {}^{n+1} \mathbf C_{\text{i}} + \tilde{\beta} \ {}^{n+1}\mathbf C_{\text{1i}},
\end{equation}
where $\tilde{\beta} \in \mathbb{R}$ is unknown, since it depends on ${}^{n+1}\mathbf C_{\text{1i}}$.
In order to enforce the incompressibility condition $\det({}^{n+1}\mathbf C_{\text{1i}})=1$,
we modify the right-hand side by introducing an additional term $\tilde{\varepsilon} \ {}^{n+1}\mathbf C_{\text{1i}}$
\begin{equation}\label{UpdateCii3}
{}^{n+1} \mathbf C_{\text{1i}} = {}^{n} \mathbf C_{\text{1i}} +
\xi \ \varkappa_1 \ c_1 \ {}^{n+1} \mathbf C_{\text{i}} + \tilde{\beta} \ {}^{n+1}\mathbf C_{\text{1i}} + \tilde{\varepsilon} \ {}^{n+1}\mathbf C_{\text{1i}},
\end{equation}
where $\tilde{\varepsilon}$ is determined using the incompressibility condition.
In analogy to the previous subsection, this additional term is chosen in such a way as not to spoil the weak invariance
of the solution. Rearranging the terms, we arrive at
\begin{equation}\label{UpdateCii4}
\tilde{z} \ {}^{n+1} \mathbf C_{\text{1i}} = {}^{n} \mathbf C_{\text{1i}} +
\xi \ \varkappa_1 \ c_1 \ {}^{n+1} \mathbf C_{\text{i}},
\end{equation}
where $\tilde{z}:=1 - \tilde{\beta} - \tilde{\varepsilon}$ is computed such that $\det({}^{n+1}\mathbf C_{\text{1i}})=1$.
The solution for this problem takes a very simple form
\begin{equation}\label{UpdateCii5}
{}^{n+1} \mathbf C_{\text{1i}} = \overline{{}^{n} \mathbf C_{\text{1i}} +
\xi \ \varkappa_1 \ c_1 \ {}^{n+1} \mathbf C_{\text{i}}}.
\end{equation}

\textbf{Remark 3.} Interestingly, this update formula has exactly the same structure as the
explicit update formula for the Maxwell fluid which was reported in \cite{ShutovLandgraf2013}.
This similarity is due to the fact that the evolution equation $\eqref{ReducedEvolut2}_1$
coincides with the evolution equation for the Maxwell fluid, in case of constant $\xi$.
In contrast to the derivation in \cite{ShutovLandgraf2013}, the derivation in the current section exploits
the notion of the weak invariance as a guideline. $\Box$

In exactly the same way, the following update formula is obtained for ${}^{n+1} \mathbf C_{\text{2i}}$:
\begin{equation}\label{UpdateCii6}
{}^{n+1} \mathbf C_{\text{2i}} = \overline{{}^{n} \mathbf C_{\text{2i}} +
\xi \ \varkappa_2 \ c_2 \ {}^{n+1} \mathbf C_{\text{i}}}.
\end{equation}
In a compact form, update formulas \eqref{UpdateCii5} and \eqref{UpdateCii6} are rewritten as
\begin{equation}\label{UpdateCii7}
{}^{n+1} \mathbf C_{\text{1i}} = \mathfrak{C}_{\text{1i}} ({}^{n+1} \mathbf C_{\text{i}}, \xi), \quad
{}^{n+1} \mathbf C_{\text{2i}} = \mathfrak{C}_{\text{2i}} ({}^{n+1} \mathbf C_{\text{i}}, \xi).
\end{equation}

\subsection{Overall procedure: partitioned Euler Backward method (PEBM)}

Multiplying $\eqref{prob6}_1$ with $\Delta t$, we arrive at
the following incremental consistency condition for finding $\xi$
\begin{equation}\label{Consistency}
\xi \eta = \displaystyle \Delta t \Big\langle \frac{\displaystyle 1}{\displaystyle f_0}
f \Big\rangle^{m}, \quad \text{where} \ f =\widetilde{f} ({}^{n+1} \mathbf C_{\text{i}}, {}^{n+1} \mathbf C_{\text{1i}},
 {}^{n+1} \mathbf C_{\text{2i}}, \xi).
\end{equation}
Here, the overstress function $\widetilde{f} ({}^{*} \mathbf C_{\text{i}}, {}^{*}
\mathbf C_{\text{1i}}, {}^{*}
\mathbf C_{\text{2i}}, {}^{*}\xi)$ is determined from $\eqref{prob6}_2$ and $\eqref{prob6}_3$ as follows
\begin{equation}\label{Overstress123}
\widetilde{f} ({}^{*} \mathbf C_{\text{i}}, {}^{*}
\mathbf C_{\text{1i}}, {}^{*}
\mathbf C_{\text{2i}}, {}^{*}\xi):=
\mathfrak{F}_1 ({}^{*} \mathbf C_{\text{i}}, {}^{*}
\mathbf C_{\text{1i}}, {}^{*}
\mathbf C_{\text{2i}}) - \sqrt{\frac{2}{3}} \big[ K +  {}^{n+1} R ({}^{*}\xi) \big],
\end{equation}
\begin{equation}\label{Overstress124}
\mathfrak{F}_1 ({}^{*} \mathbf C_{\text{i}}, {}^{*}
\mathbf C_{\text{1i}}, {}^{*}
\mathbf C_{\text{2i}}):= \sqrt{\text{tr}
\big[ \big( {}^{n+1} \mathbf C \ {}^{*}\tilde{\mathbf T} -
{}^{*}\mathbf C_{\text{i}} \ {}^{*}\tilde{\mathbf X} \big)^{2} \big] },
\end{equation}
\begin{equation}\label{Overstress125}
 {}^{*}\tilde{\mathbf T} ({}^{*} \mathbf C_{\text{i}}) :=  \mu \ {}^{n+1}\mathbf C^{-1} (\overline{{}^{n+1} \mathbf C} \
{}^{*}\mathbf C_{\text{i}}^{-1})^{\text{D}},
\end{equation}
\begin{equation}\label{Overstress126}
 {}^{*} \tilde{\mathbf X} ({}^{*} \mathbf C_{\text{i}}, {}^{*}
\mathbf C_{\text{1i}}, {}^{*}
\mathbf C_{\text{2i}}) := \frac{c_1}{2} \  {}^{*}\mathbf C_{\text{i}}^{-1} \
( {}^{*}\mathbf C_{\text{i}} \  {}^{*}\mathbf C_{\text{1i}}^{-1})^{\text{D}} +
\frac{c_2}{2} \  {}^{*}\mathbf C_{\text{i}}^{-1} \
( {}^{*}\mathbf C_{\text{i}} \  {}^{*}\mathbf C_{\text{2i}}^{-1})^{\text{D}}.
\end{equation}

The following predictor-corrector scheme is implemented in this study:
\begin{itemize}
\item[1]   \textbf{Elastic predictor:} Evaluate the trial overstress for frozen internal variables
\begin{equation}\label{OverstressTrial}
 {}^{\text{trial}} f:= \widetilde{f} ({}^{n} \mathbf C_{\text{i}}, {}^{n} \mathbf C_{\text{1i}},
 {}^{n} \mathbf C_{\text{2i}}, 0).
\end{equation}
If $ {}^{\text{trial}} f \leq 0$ then the current stress state lies in the elastic domain. Put $\xi =0$,
${}^{n+1} \mathbf C_{\text{i}} = {}^{n} \mathbf C_{\text{i}}$,
${}^{n+1} \mathbf C_{\text{1i}} = {}^{n} \mathbf C_{\text{1i}}$,
${}^{n+1} \mathbf C_{\text{2i}} = {}^{n} \mathbf C_{\text{2i}}$, ${}^{n+1} s = {}^{n} s$, ${}^{n+1} s_{\text{d}} = {}^{n} s_{\text{d}}$.
The time step is thus complete.
If $ {}^{\text{trial}} f > 0$ then proceed to the plastic corrector step.
\item[2]  \textbf{Plastic corrector:} The corrector step consists of the following substeps.

\begin{itemize}

\item[2.1] The initial (rough) estimation for ${}^{n+1}  \mathbf{C}_{\text{i}}$ is obtained from  ${}^{n}  \mathbf{C}_{\text{i}}$
by the same shift (push-forward operation) which brings $\overline{{}^{n}  \mathbf{C}}$ to $\overline{{}^{n+1}  \mathbf{C}}$. More precisely:
\begin{equation}\label{DefineShift}
{}^{\text{est}}  \mathbf{C}_{\text{i}} \gets \mathbf{F}^{-\text{T}}_{\text{sh}} \ {}^{n}  \mathbf{C}_{\text{i}} \ \mathbf{F}^{-1}_{\text{sh}}, \ \text{where} \
\mathbf{F}_{\text{sh}} := (\overline{{}^{n+1}  \mathbf{C}^{-1} \  {}^{n}  \mathbf{C}})^{1/2}.
\end{equation}
Next, estimate $\xi >0$ by solving the following scalar equation
\begin{equation}\label{RoughXi}
{}^{\text{est}}\xi =  ({}^{\text{trial}} f  - ({}^{\text{est}}\xi \eta/ \Delta t)^{1/m} )/(2 \mu).
\end{equation}

\item[2.2]
The internal variables ${}^{\text{n+1}} \mathbf{C}_{\text{1i}}$ and ${}^{\text{n+1}} \mathbf{C}_{\text{2i}}$
are estimated using the explicit update formulas
\begin{equation}\label{RoughCii}
{}^{\text{est}} \mathbf{C}_{\text{1i}} \gets \mathfrak{C}_{\text{1i}} ({}^{\text{est}} \mathbf C_{\text{i}}, {}^{\text{est}}\xi), \quad
{}^{\text{est}} \mathbf{C}_{\text{2i}} \gets \mathfrak{C}_{\text{2i}} ({}^{\text{est}} \mathbf C_{\text{i}}, {}^{\text{est}}\xi).
\end{equation}

\item[2.3] Resolve the following incremental consistency condition with respect to $\xi$
\begin{equation}\label{Consistency2}
\xi \eta = \displaystyle \Delta t \Big\langle \frac{\displaystyle 1}{\displaystyle f_0}
\widetilde{f} (\mathfrak{C}_{\text{i}}({}^{\text{est}}\mathbf C_{\text{1i}}, {}^{\text{est}}\mathbf C_{\text{2i}}, \xi),
{}^{\text{est}} \mathbf C_{\text{1i}}, {}^{\text{est}} \mathbf C_{\text{2i}}, \xi) \Big\rangle^{m}.
\end{equation}
In other words, a time step is performed with fixed $\mathbf C_{\text{1i}} \equiv {}^{\text{est}}\mathbf C_{\text{1i}}$ and
$\mathbf C_{\text{2i}} \equiv {}^{\text{est}}\mathbf C_{\text{2i}}$.
Let $ \tilde{\xi}$ be the solution of \eqref{Consistency2}. Using it, compute
\begin{equation}\label{Substep22}
{}^{\text{est}} \xi \gets \tilde{\xi}, \quad
{}^{\text{est}}\mathbf C_{\text{i}} \gets \mathfrak{C}_{\text{i}}({}^{\text{est}}\mathbf C_{\text{1i}},
{}^{\text{est}}\mathbf C_{\text{2i}},  \tilde{\xi}).
\end{equation}

\item[2.4]

For better accuracy, repeat substeps 2.2, 2.3, 2.2, 2.3; then assign
\begin{equation}\label{Substep24}
\xi \gets {}^{\text{est}} \xi, \quad {}^{n+1}  \mathbf C_{\text{i}} \gets {}^{\text{est}}  \mathbf C_{\text{i}},
\quad {}^{n+1}  \mathbf C_{\text{1i}} \gets {}^{\text{est}}  \mathbf C_{\text{1i}}, \quad
\quad {}^{n+1}  \mathbf C_{\text{2i}} \gets {}^{\text{est}}  \mathbf C_{\text{2i}}.
\end{equation}
Finally, the variables $s$ an $s_{\text{d}}$ are updated by  \eqref{updateSSd}.
The plastic corrector step is thus complete.
\end{itemize}
\end{itemize}

\textbf{Remark 4.} In the current algorithm, substeps 2.2 and 2.3 are repeated for three times.
In the general context of the decoupled Euler Backward method,
these repetitions correspond to the so-called relaxation iterations.
As an alternative to the proposed algorithm, one may consider a full relaxation
which consists in numerous repetitions of 2.2 and 2.3, carried to convergence.
As will be shown in the following,
the full relaxation is not necessary and \emph{only three relaxation iterations
yield very good results} for different metallic materials.
Such accurate results are achieved due to the preliminary substep 2.1,
which provides a good approximation of the solution
even for finite values of $\xi$.
Without the preliminary push-forward \eqref{DefineShift},
the integration procedure would be less accurate. $\Box$

Instead of solving a system of nonlinear algebraic equations with
respect to ${}^{n+1} \mathbf C_{\text{i}}$, ${}^{n+1} \mathbf C_{\text{1i}}$, ${}^{n+1} \mathbf C_{\text{2i}}$, and $\xi$,
as it was carried out in previous studies, only a scalar
consistency equation has to be solved with respect to $\xi$. Namely, equation
\eqref{Consistency2} has to be resolved for three times.
The resulting numerical scheme is first order accurate\footnote{The
reader interested in higher-order integration methods for elasto-plasticity
is referred to \cite{HartmannBier2008, Eidel, EidelStumpf} and references cited therein.}; the solution remains
bounded even for very large time steps.
Moreover, the geometric property \eqref{geopr0} is exactly satisfied and the positive definiteness
of $\mathbf C_{\text{i}}$ and  $\mathbf C_{\text{ii}}$ is
guaranteed even for very large time steps and strain increments.

The push-forward operation \eqref{DefineShift} exactly retains the weak invariance of the solution (see Appendix D).
For that reason, the numerical solution exhibits the same weak invariance property as the original continuum model.

\section{Simulation results}

\subsection{Constitutive parameters for real materials}

In this section, the novel PEBM is tested using two different sets of constitutive parameters
pertaining to two different materials.
The accuracy of the stress computation by PEBM will be compared with the
Euler Backward method with subsequent correction (EBMSC) and the exponential method (EM). These
conventional methods are shortly summarized in Appendix E.

\subsubsection{Material parameters for 5754-O aluminum alloy}

\begin{figure}\centering
\scalebox{0.95}{\includegraphics{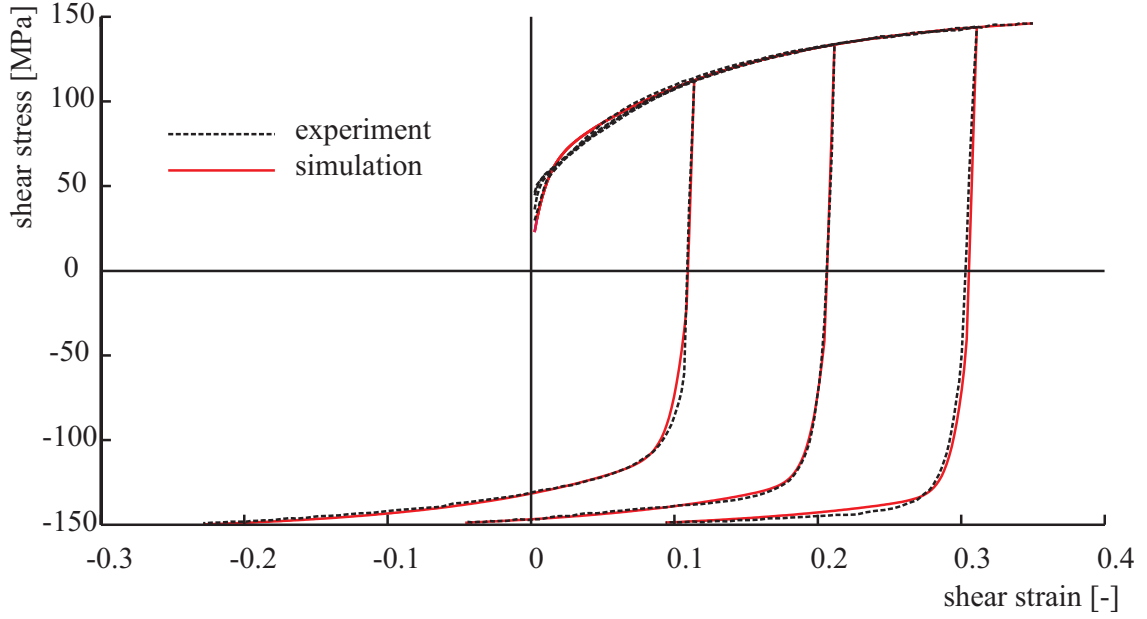}} \caption{Experimental and theoretical results for the stress response of 5754-O aluminum alloy
subjected to simple shear.
Experimental data reported in \cite{Chaparro2008, LaurentGreeze2009}; simulations performed using parameters from table \ref{tab1}.
\label{figFor5754}}
\end{figure}

We adopt the experimental data reported for 5754-O aluminum alloy in \cite{Chaparro2008, LaurentGreeze2009}.
Four different shear tests are employed for the parameter identification: three Bauschinger tests with different prestrain levels
and a monotonic test. For this material, the rate-dependence is weak \cite{LaurentGreeze2009}; therefore, we switch-off
the viscosity by putting $\eta =0$. The identified constitutive parameters are summarized in table \ref{tab1}. The comparison of the
theoretical and experimental results is shown in Figure \ref{figFor5754}. As can be seen from the figure, the model
with two parameters for the isotropic hardening and four parameters for the kinematic hardening describes the real stress response with
a good fidelity.

\begin{table}[h!]
\caption{Set of constitutive parameters for 5754-O aluminum alloy}
\begin{center}
\begin{tabular}{| l l l l l |}
\hline
$k$ [MPa] & $\mu$ [MPa]  &  $c_1$  [MPa]  & $c_2$  [MPa]  &  $\gamma$ [MPa] \\ \hline
68630     & 26310        &   115.5        & 11500         & 1963            \\ \hline
\end{tabular} \\
\begin{tabular}{|l l l l l l|}
\hline
$K$ [MPa] & $m$ [-] & $\eta$ [$\text{s}$]    & $\varkappa_1$ [$\text{MPa}^{-1}$] & $\varkappa_2$ [$\text{MPa}^{-1}$] & $\beta$ [-]   \\ \hline
31.5      & 1.0     & 0                      &  0.0885                           &  0.01676                           &  13.33        \\ \hline
\end{tabular}
\end{center}
\label{tab1}
\end{table}

\subsubsection{Material parameters for 42CrMo4 steel}

\begin{figure}\centering
\scalebox{0.95}{\includegraphics{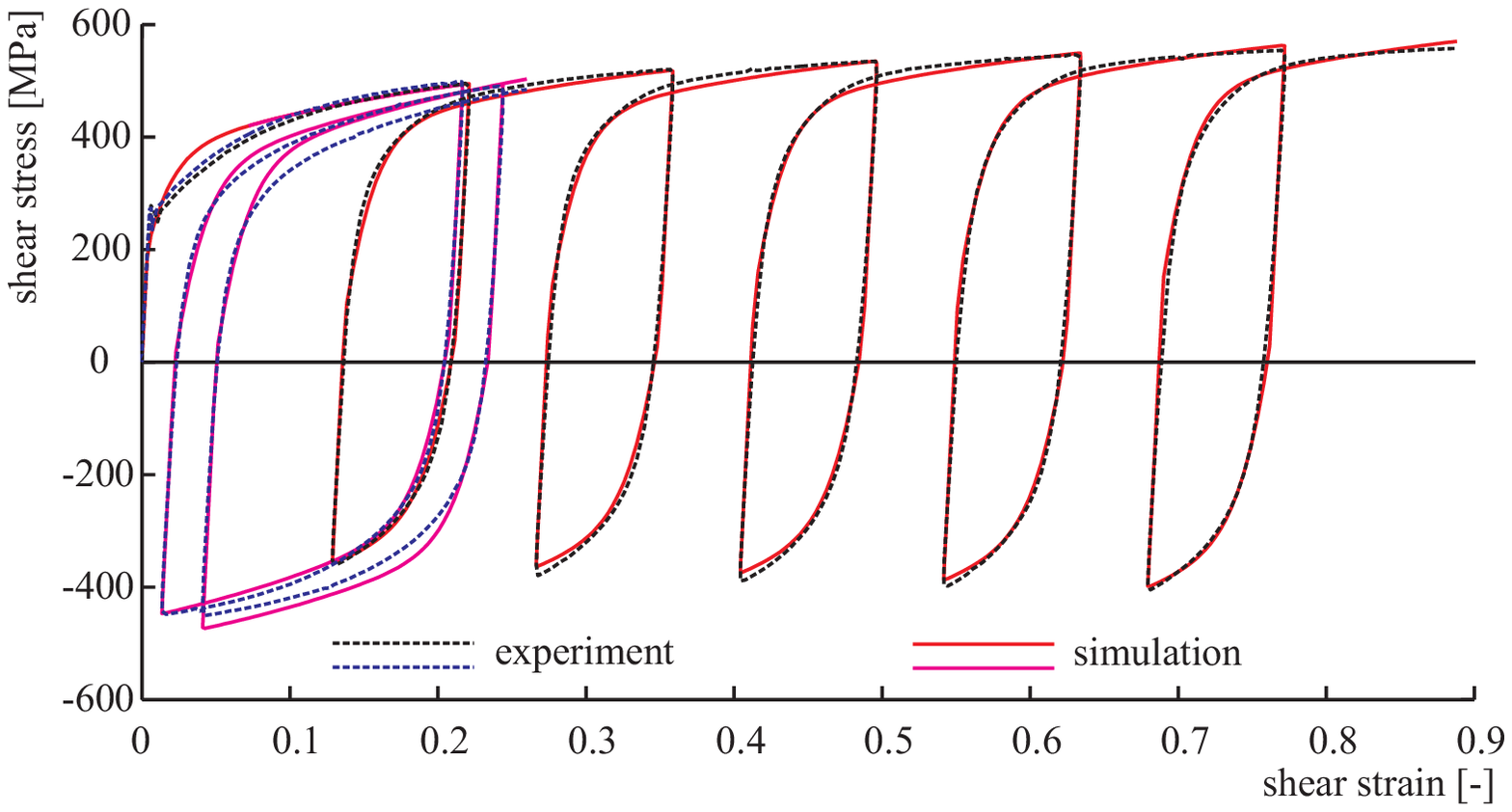}} \caption{Experimental and theoretical results for the stress response of 42CrMo4 steel;
two cyclic torsion tests are presented.
Experimental data reported in \cite{ShutovKuprin}; simulations performed using parameters from table \ref{tab2}.
\label{figFor42CrMo4}}
\end{figure}

The experimental data for 42CrMo4 steel were reported in \cite{ShutovKuprin}. For details concerning
the experimental set up and simulation of the torsion test, the reader is referred to \cite{ShutovKuprin}. All
 tests are performed with a constant shear
strain rate $| \dot{\gamma}_{\text{shear}} | = 0.07 s^{-1}$.
Two different cyclic torsion tests  are adopted in the current study for the identification of
material parameters. The identified parameters are summarized in table \ref{tab2}. As can be seen from Figure \ref{figFor42CrMo4}, the considered model
describes the nonlinear stress response upon the load reversal with a good accuracy.

\begin{table}[h!]
\caption{Set of constitutive parameters for 42CrMo4 steel}
\begin{center}
\begin{tabular}{| l l l l l |}
\hline
$k$ [MPa] & $\mu$ [MPa]  &  $c_1$  [MPa]  & $c_2$  [MPa]  &  $\gamma$ [MPa] \\ \hline
135600     & 52000        &   1674         & 22960         &   145           \\ \hline
\end{tabular} \\
\begin{tabular}{|l l l l l l|}
\hline
$K$ [MPa] & $m$ [-] & $\eta$ [$\text{s}$]    & $\varkappa_1$ [$\text{MPa}^{-1}$] & $\varkappa_2$ [$\text{MPa}^{-1}$] & $\beta$ [-]   \\ \hline
335       & 2.26    & 500 000                 &  0.00386                           & 0.0043102                         &   0.0503       \\ \hline
\end{tabular}
\end{center}
\label{tab2}
\end{table}

For all the following tests we switch-off the viscous effects: material parameters from table \ref{tab2} are taken for 42CrMo4 steel
with $\eta=0$ and $m=1$.

We note that the material parameters which appear in tables \ref{tab1} and \ref{tab2}
correspond to a very fast saturation of nonlinear kinematic hardening with large backstresses in the saturated state.
From the numerical standpoint, such a type of stress response is highly unfavorable:
much larger computational errors are expected for these parameters than in case of slow saturation with
small backstresses. Thus, the following numerical tests are assumed to be representative for a broad class of metallic materials.

\subsection{Accuracy test: non-proportional loading}

Consider the time instancies $T_0=0$s, $T_1=100$s, $T_2=200$s, $T_3=300$s.
Let us simulate the stress response
in the time interval $t \in [T_0,T_3]$
for the following local loading program
\begin{equation}\label{loaprog0}
\mathbf F (t) = \overline{\mathbf F^{\prime} (t)},
\end{equation}
where the auxiliary function $\mathbf F^{\prime}(t)$ is given by
\begin{equation*}\label{loaprog}
\mathbf F^{\prime} (t) :=
\begin{cases}
    (1 - t/T_1) \mathbf F_1  + (t/T_1) \mathbf F_2 \quad \text{if} \ t \in [T_0,T_1] \\
    (2 - t/T_1) \mathbf F_2  + (t/T_1-1) \mathbf F_3 \quad \text{if} \ t \in (T_1,T_2] \\
     (3 - t/T_1) \mathbf F_3  + (t/T_1-2) \mathbf F_4 \quad \text{if} \ t \in (T_2,T_3]
\end{cases}.
\end{equation*}
Here, $\mathbf F_1$, $\mathbf F_2$, $\mathbf F_3$, and $\mathbf F_4$ are the key points, which are defined in a fixed
Cartesian coordinate system through
\begin{equation*}\label{loaprog2}
\mathbf F_1 :=\mathbf 1, \
\mathbf F_2 := \left(
\begin{array}{ccc}
2 & 0 & 0 \\
0 & \displaystyle \frac{1}{\sqrt2} & 0 \\
0 & 0 & \displaystyle \frac{1}{\sqrt2}
\end{array}
\right), \
\mathbf F_3 := \left(
\begin{array}{ccc}
1 & 1  & 0 \\
0 & 1 & 0 \\
0 & 0 & 1
\end{array}
\right), \
\mathbf F_4 := \left(
\begin{array}{ccc}
\displaystyle \frac{1}{\sqrt2} & 0 & 0 \\
0 & 2 & 0 \\
0 & 0 & \displaystyle \frac{1}{\sqrt2}
\end{array}
\right).
\end{equation*}
This program defines a rather arbitrary loading, representative for metal forming applications.
Note that this loading program exhibits an abrupt change of the load path at $t=100 \text{s}$ and $t=200 \text{s}$.
For the numerical test we impose the following initial conditions, corresponding to unloaded isotropic initial state:
\begin{equation}\label{Inico}
\mathbf C_{\text{i}}|_{t=0} =  \mathbf 1, \quad    \mathbf C_{\text{1i}}|_{t=0}=  \mathbf 1,
\quad  \mathbf C_{\text{2i}}|_{t=0}=  \mathbf 1, \quad s|_{t=0}=0, \quad s_{\text{d}}|_{t=0}=0.
\end{equation}
The numerical solution obtained with extremely small time step ($\Delta t = 0.0025 \text{s}$) using the
EBMSC will be seen as the exact solution.
The numerical error function is computed as a distance between the exact and the numerically computed Cauchy stresses
\begin{equation}\label{ErrorDef}
\text{Error}(t): = \| \mathbf T^{\text{numerical}}(t) - \mathbf T^{\text{exact}}(t)  \|.
\end{equation}

\begin{figure}\centering
\scalebox{0.8}{\includegraphics{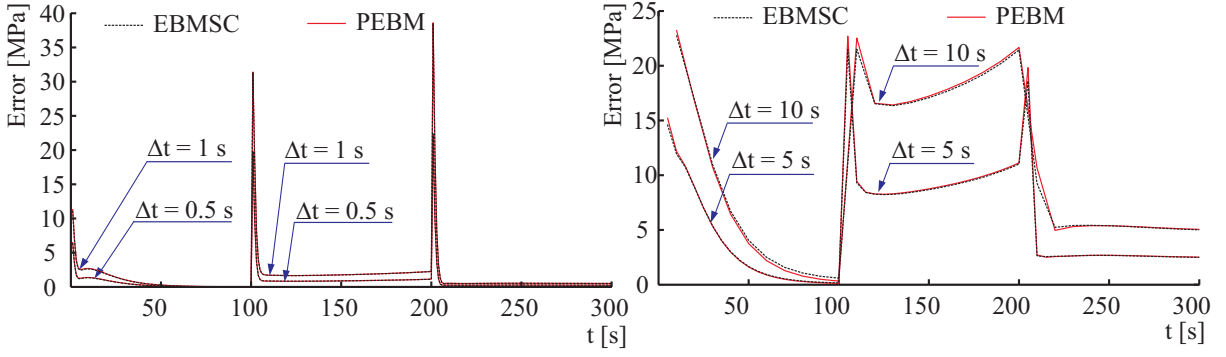}} \caption{Dependence of the integration error on the time step size $\Delta t$.
Material parameters of 5754-O aluminum alloy are used.
\label{ErrorForAl}}
\end{figure}

\begin{figure}\centering
\scalebox{0.8}{\includegraphics{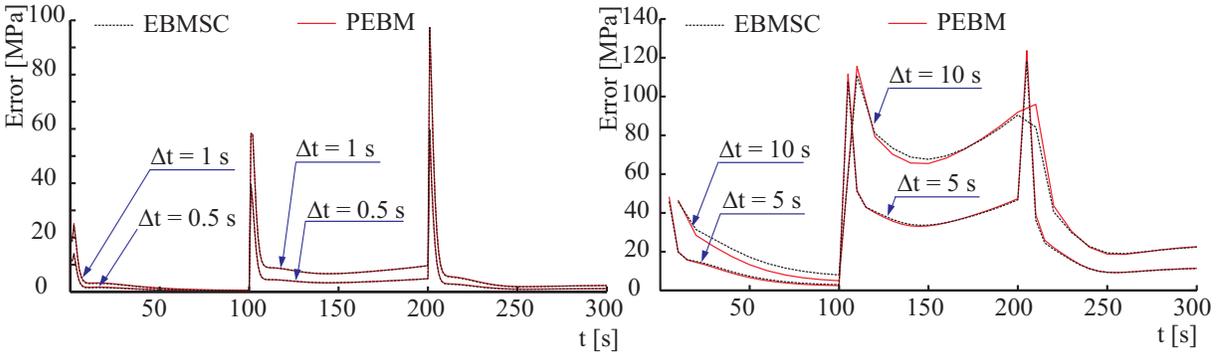}} \caption{Dependence of the integration error on the time step size $\Delta t$.
Material parameters of 42CrMo4 steel are used.
\label{ErrorForSt}}
\end{figure}

The error function is plotted for different values of the time step size $\Delta t$;
Figures \ref{ErrorForAl} and \ref{ErrorForSt} correspond to material parameters for Al 5754O and 42CrMo4 steel, respectively.
The numerical tests reveal that the proposed partitioned Euler Backward method (PEBM) and the
conventional Euler Backward with subsequent correction (EBMSC)
have a similar integration error even
for very big time steps $\Delta t = 10 \ \text{s}$. Both algorithms are first order
accurate and they are \emph{comparable in accuracy}, but
the computational effort for the new algorithm is much smaller than for EBMSC. A more detailed discussion
of the computational costs will be presented in Section 4.4.

\subsection{Accuracy test: iso-error maps}

The iso-error maps are useful for the visualization of the integration error within a single time step \cite{SimoHughes, LeeKim2015}.
In this subsection we consider only incompressible deformation paths which can be described by
\begin{equation}\label{DefPathIsoError}
\mathbf F (t) = F_{11}(t) \ \mathbf e_1 \otimes \mathbf e_1 + (F_{11}(t))^{-1/2}
(\mathbf e_2 \otimes \mathbf e_2 + \mathbf e_3 \otimes \mathbf e_3) +
F_{12}(t) \ \mathbf e_1 \otimes \mathbf e_2.
\end{equation}
Just as in the previous subsection, the viscous effects are switched off.
Thus, the stress response is rate independent. For simplicity, the time $t$ is understood
here as a non-dimensional monotonic loading parameter.
Moreover, we switch off the isotropic hardening as well, by putting $\gamma = 0$ for both materials.
To mimic the effect of isotropic hardening, we put $K = 178.8$ MPa for the AA 5754-O and $K = 400$ MPa for the 42CrMo4 steel.

\begin{figure}\centering
\scalebox{0.85}{\includegraphics{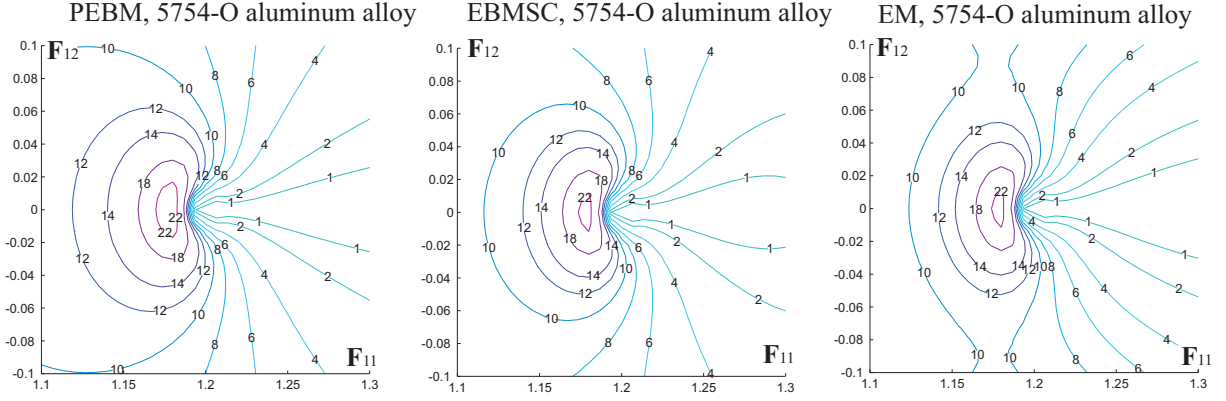}} \caption{Iso-error maps for different integration methods using material parameters
of 5754-O aluminum alloy. Prestrain in tension.
\label{IsoErrorTensionAlu}}
\end{figure}

\begin{figure}\centering
\scalebox{0.85}{\includegraphics{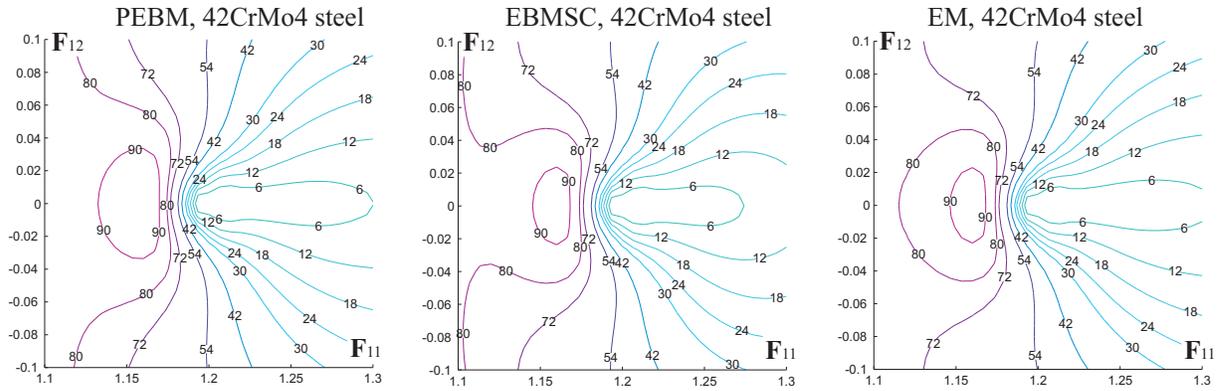}} \caption{Iso-error maps for different integration methods using material parameters
of 42CrMo4 steel. Prestrain in tension.
\label{IsoErrorTensionSteel}}
\end{figure}

First, we consider a 20\% prestrain in uniaxial tension
\begin{equation}\label{UniaxTensPrestrain}
F_{11}(t) = 1 + 0.2 t, \quad F_{12}(t)=0, \quad t \in [0,1].
\end{equation}
This prestrain is simulated by 500 time steps.
After that, a single time step is performed to a new state, characterized by $(F_{11},F_{12})$.
The numerical solution obtained for this step
using EBMSC with 300 subincrements is considered as the exact solution.
Using this exact solution, the error is defined as in the previous subsection (see equation \eqref{ErrorDef}).
The iso-error maps for the new algorithm (PEBM), for the Euler Backward method with subsequent correction (EBMSC),
and for the exponential map (EM) are presented in Figure \ref{IsoErrorTensionAlu} for the  5754-O aluminum alloy
and in Figure \ref{IsoErrorTensionSteel} for the 42CrMo4 steel. As in the previous subsection,
the computational methods exhibit a similar numerical error. For all three methods, the error is relatively small in the recent loading direction
and large in the opposite direction. The PEBM-error is slightly larger than the EBMSC-error for strain increments opposite
to the recent loading direction and slightly smaller for increments in the loading direction.

\begin{figure}\centering
\scalebox{0.85}{\includegraphics{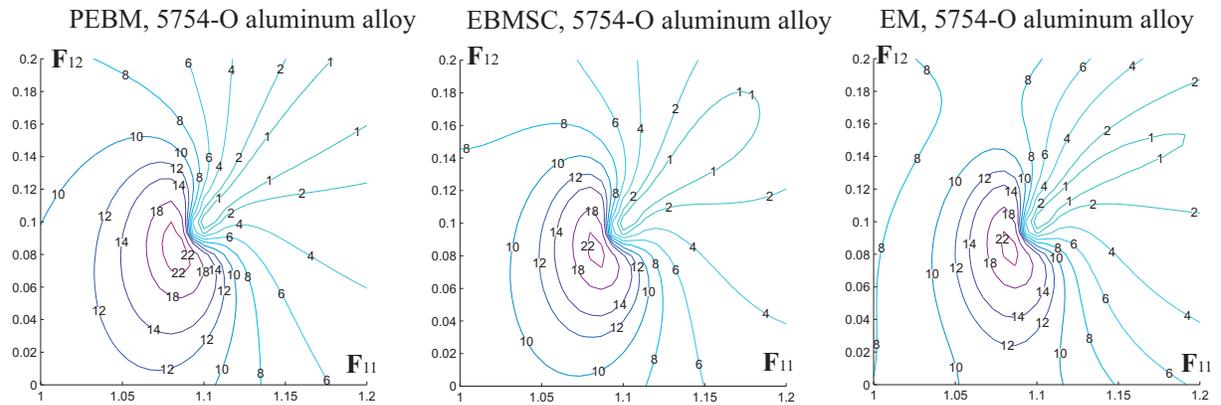}} \caption{Iso-error maps for different integration methods using material parameters
of 5754-O aluminum alloy. Combined prestrain in tension and shear.
\label{IsoErrorTensionAndShearAlu}}
\end{figure}

\begin{figure}\centering
\scalebox{0.85}{\includegraphics{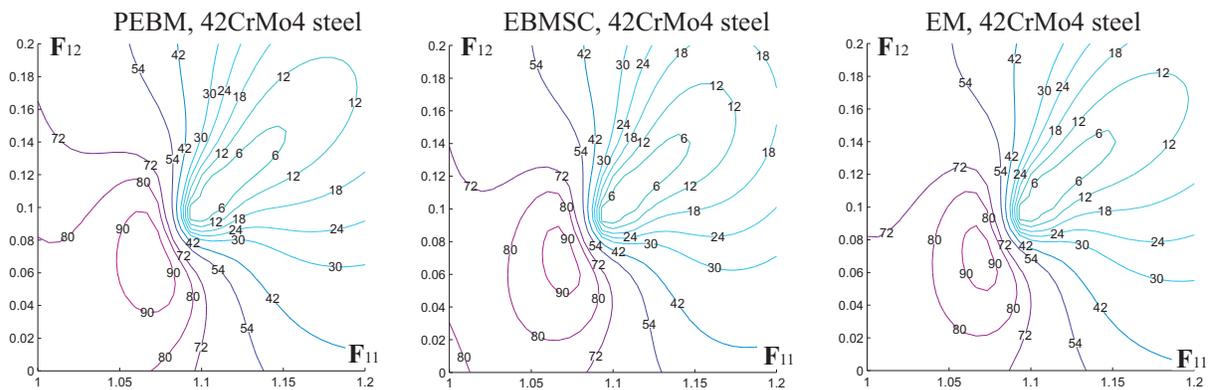}} \caption{Iso-error maps for different integration methods using material parameters
of 42CrMo4 steel. Combined prestrain in tension and shear.
\label{IsoErrorTensionAndShearSteel}}
\end{figure}

Next, we consider a combined prestrain with 10\% tension and 10\% shear
\begin{equation}\label{UniaxTensPrestrain}
F_{11}(t) = 1 + 0.1 t, \quad F_{12}(t)=0.1 t, \quad t \in [0,1].
\end{equation}
The iso-error maps for free different integration algorithms are shown in Figure \ref{IsoErrorTensionAndShearAlu} for the  5754-O aluminum alloy
and in Figure \ref{IsoErrorTensionAndShearSteel} for the 42CrMo4 steel. Just as in case of uniaxial prestrain, all three algorithms are equivalent regarding
the accuracy of the stress computation. Again, small errors occur in the recent loading direction along with relatively large errors in the opposite direction.

\subsection{Analysis of computational efforts}

For the novel PEBM, the computational effort depends solely on the number of Newton iterations required to resolve the
incremental consistency condition \eqref{Consistency2}.
For the conventional EBMSC and EM, a system of 19 coupled nonlinear
equations has to be resolved at each time step (cf. Appendix E).
The simultaneous numerical solution
of these equations is problematic for large strain increments. In order to obtain
a more reliable procedure, we employ the return mapping approach.
First, a numerical procedure is created which computes
the unknown
${}^{n+1} \mathbf C_{\text{i}}$, ${}^{n+1} \mathbf C_{\text{1i}}$, and ${}^{n+1} \mathbf C_{\text{2i}}$
as functions of $\xi$. We call this a ``$\xi$-step".
Each $\xi$-step is based on the iterative process of the Newton-Raphson type where 18 nonlinear equations have to be solved.
Then, the correct value of $\xi$ is found as a solution of the incremental consistency condition \eqref{AppendixE}. Just as for PEBM, this
scalar equation is solved using the Newton method. For the novel PEBM each Newton
iteration requires a finite number of operations on 3x3 matrices (cf. Sections 3.2 and 3.3). These operations include the matrix
multiplication, computation of the square root and computation of the inverse.
On the other hand, for the conventional EBMSC and EM each Newton iteration calls the $\xi$-step, which is much more expensive since
an 18x18 matrix has to be built and inverted several times within the iterative process. Generally, the Newton iteration for the PEBM
is approximately two orders of magnitude cheaper than the Newton iteration for the EBMSC or EM. 

Since the
$\xi$-step involves the iterative process, it is (a priory) less reliable than a closed-form solution used in PEBM.
For $\xi > 0.1$, the straightforward implementation of the Newton-Raphson method to the $\xi$-step yields a divergent
iteration process. In order to obtain a reliable numerical procedure
for EBMSC and EM, the $[0, \xi]$ interval is subincremented, which further increases the computational costs.

\begin{figure}\centering
\scalebox{0.85}{\includegraphics{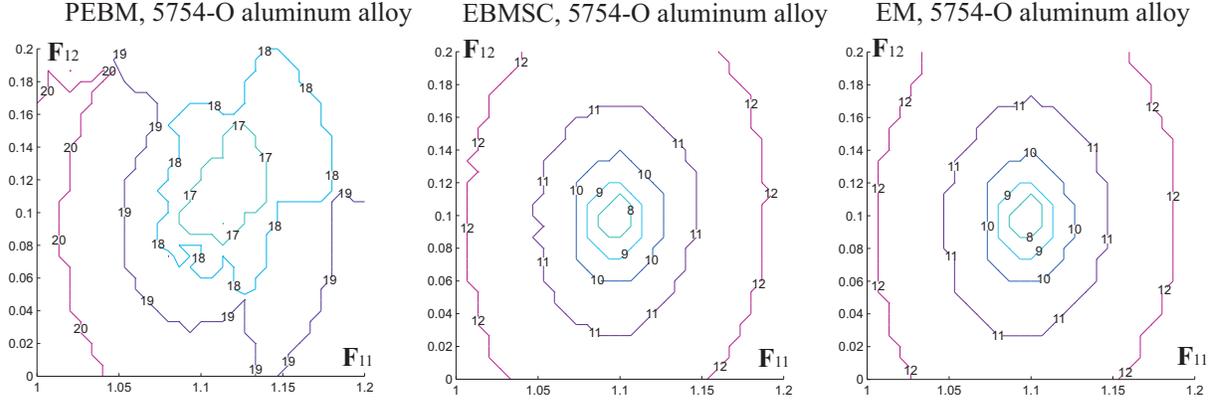}}
\caption{Isolines for the number of Newton iterations required by different methods.
\label{figIterNumber}}
\end{figure}

The setting used in the previous subsection can be also adopted for the analysis of the computational efforts.
Toward that end we consider a combined prestrain with 10\% tension and 10\% shear, as described by \eqref{UniaxTensPrestrain}.
The prestrain is simulated numerically by 500 time steps. After that, a single time step is performed to a new $(F_{11},F_{12})$-state.
The number of Newton iterations required by PEBM, EBMSC and EM to make this step is plotted in Fig. \ref{figIterNumber}.
Recall that for the PEBM the incremental consistency condition has to be solved for three times,
since three relaxation iterations are implemented in the current study (cf. Remark 4). Within the first relaxation iteration, the initial
approximation for the unknown $\xi$ is zero. For the subsequent relaxation iterations, the solution from the previous relaxation iterations acts
as the initial approximation. The total (accumulated) number of Newton iterations is plotted for the PEBM.
Although PEBM requires more iterations than the conventional methods, the difference in the iteration number is not large. Therefore,
PEBM is still much more efficient than EBMSC or EM. In contrast to the conventional EBMSC and EM, \emph{no convergence issues} 
are observed for the novel PEBM.

\subsection{FEM solutions of boundary value problems}

The proposed PEBM is implemented into the commercial FEM code MSC.MARC employing the Hypela2 interface.
In order to demonstrate the correctness of the numerical implementation, we
consider two different bending-dominated problems with finite displacements and rotations.
In both problems, the material parameters correspond to the 5754-O aluminum alloy (cf. table \ref{tab1}).

\begin{figure}\centering
\scalebox{0.9}{\includegraphics{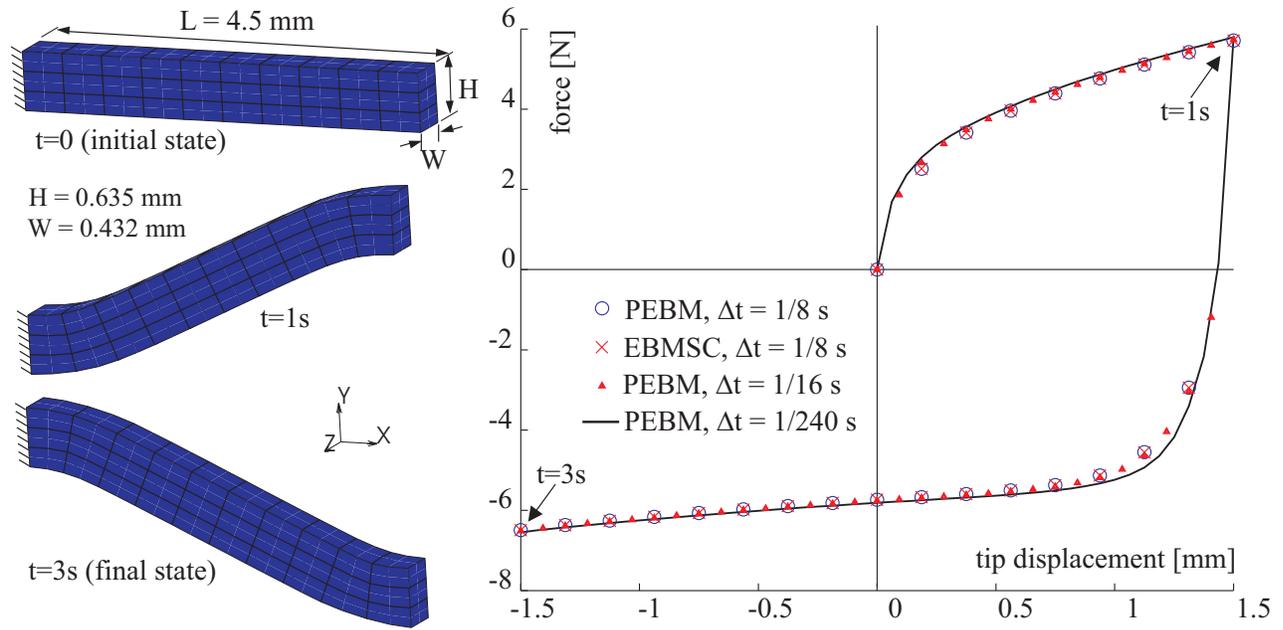}} \caption{FEM solution in MSC.MARC: displacement-controlled loading of a cantilever beam.
\label{BeamMarc}}
\end{figure}

First, a displacement-controlled problem is solved. A cantilever beam is modeled as shown in Figure \ref{BeamMarc} (left).
The left end of the beam is clamped; the right end is subjected
to a prescribed vertical displacement with zero rotation.
A total number of 30 elements of type Hex20 with the full integration is used.
A series of tests is performed with a constant time step size within each test;
the corresponding reaction-displacement curves are shown in Figure \ref{BeamMarc} (right).
The numerical solution obtained using 720 time steps (with  $\Delta t = 1/240$ s) serves as a reference.
As can be seen, even for the large time steps ($\Delta t = 1/8$ s, 24 times steps), the conventional EBMSC and the novel PEBM
produce almost identical results. Both methods \emph{exhibit quadratic convergence} and they require
approximately the same number of iterations.

\begin{figure}\centering
\scalebox{0.9}{\includegraphics{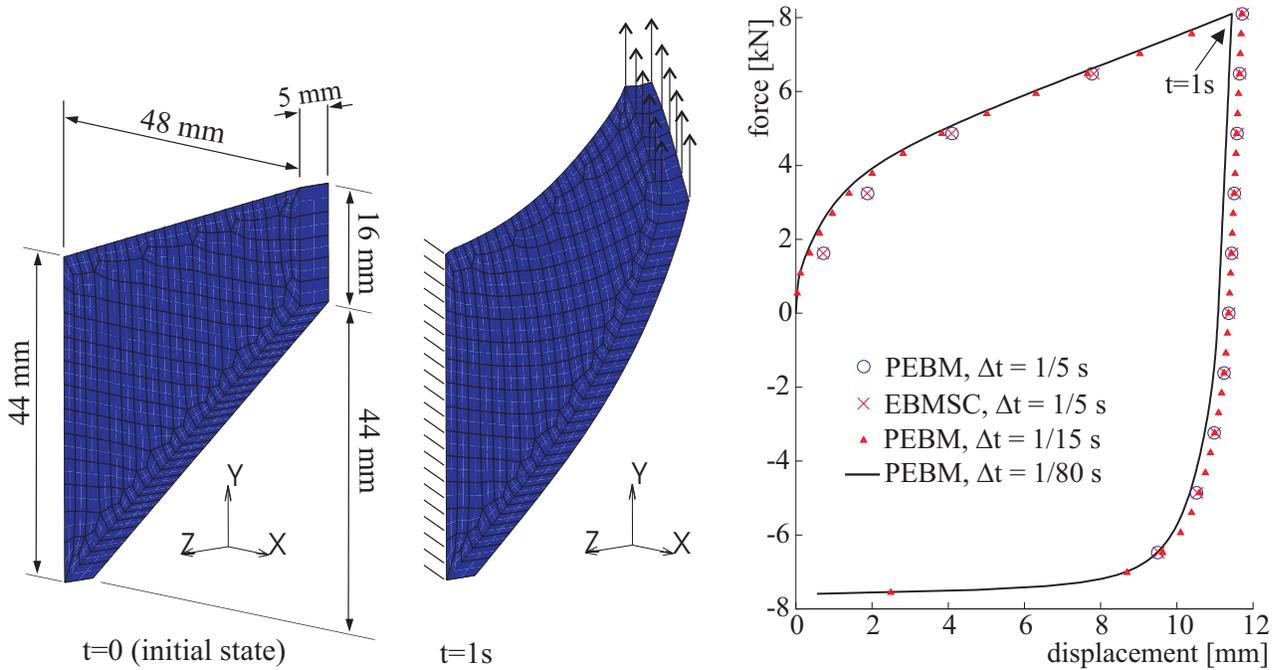}} \caption{FEM solution in
MSC.MARC: force-controlled loading of Cook's membrane.
\label{CookMarc}}
\end{figure}

Further, a non-monotonic force-controlled loading is applied to the Cook membrane as shown in Fig. \ref{CookMarc}.
The mesh contains 213 elements of type Hex20 with the full integration.
The obtained force-displacement curves are shown in Figure \ref{CookMarc} (right) for different
values of the time step size.
The numerical solution obtained using 240 time steps (with  $\Delta t = 1/80$ s) serves as a reference.
Just as in the previous problem, EBMSC and PEBM
produce almost identical results, even for the large time steps ($\Delta t = 1/5$ s). Again,
both methods exhibit quadratic convergence and they require
approximately the same number of iterations.

Concluding this subsection, the conventional EBMSC and the novel EBM seem to be equivalent regarding the accuracy and stability of
the global FE procedure, even if large strain increments are involved.

\section{Discussion and conclusion}

A new structure-exploiting numerical scheme is proposed for the implicit time stepping. 
The analyzed model of finite-strain viscoplasticity is
based on the nested multiplicative split combined with hyperelastic relations.
The main ingredients of the new algorithm are:
\begin{itemize}
\item Elastic predictor --- plastic corrector approach, also known as return mapping;
\item Use of new analytical (closed-form) solutions for elementary decoupled problems (see Section 3.2);
\item Approximation of the inelastic strain within each time step
by a special push-forward operation (see Section 3.3)
\end{itemize}

The new algorithm exactly preserves the symmetry of the tensor-valued internal variables $\mathbf{C}_{\text{i}}$,
$\mathbf{C}_{\text{1i}}$, and $\mathbf{C}_{\text{2i}}$. Moreover, these quantities remain positive definite, which
is important, since they represent some inelastic metric tensors.
The incompressibility restriction is also satisfied, which is necessary for prevention of the
error accumulation. Another important property of the proposed algorithm is the exact retention of the weak invariance
under the isochoric change of the reference configuration. \emph{It is impressive
that the retention of the weak invariance may serve as a guide in constructing efficient numerical algorithms.} For instance,
the update formulas in Sections 3.2 and 3.3 as well as the push-forward operation \eqref{DefineShift} are obtained
employing the weak invariance as a hint. The proposed PEBM is first-order accurate; the solution remains
bounded even for very large time steps.

The performance of the new algorithm is examined using the constitutive parameters of two real materials, which is a tough test because
of large values of $c = c_1 + c_2$. The larger $c$ is, the stiffer is the problem, which makes an accurate numerical
integration even more challenging. Moreover, some convergence problems are observed on the Gauss point level while using
the Euler Backward method and exponential method, if working with big
inelastic increments ($\xi > 0.1$) and large values of $c$.
\emph{At the same time, no computational difficulties are encountered within the
proposed PEBM, since the problem-causing nonlinear algebraic equations
are now solved analytically.}
The numerical tests reveal that the proposed PEBM
exhibits a similar error as the EBMSC. Both methods
produce numerical results close to the EM.
Regarding the computational efficiency, the new local algorithm
is superior to the EBMSC and EM.
For example, for a model with two back-stress tensors considered in the paper,
the EBMSC and EM require a solution of a system of 19 nonlinear equations whereas
the new PEBM allows us to reduce the number of
equations from 19 to 1. The algorithm is implemented into MSC.MARC and
two boundary-value problems are solved numerically to demonstrate the correctness
of the FE implementation. According to the test results, the conventional EBMSC and the novel PEBM
are equivalent concerning the accuracy and reliability of the global FE procedure.

Although the advocated time-stepping method is demonstrated here using a concrete model of finite strain viscoplasticity, it
can be applied to similar models with the nested multiplicative split, developed for thermoplasticity \cite{ShutovThermal}, shape memory
alloys \cite{Helm1, Helm3, ReeseChrist, FrischkornReese}, and ratchetting \cite{Feigen, ZHU2013}.
The application of the proposed method to combined kinematic/distortional hardening with a radial flow rule
reported in \cite{ShutovPaKr} is straightforward as well.
In conclusion, the presented method leads to enhanced
computational efficiency of the local
implicit time-stepping for the anisotropic finite strain plasticity, based on Lion's approach.

\section*{Acknowledgement}

This research was supported by the Russian Science Foundation (project number 15-11-20013).

\section*{Appendix A. Weak invariance of discrete evolution equations}

\subsection*{Violation of the weak invariance for $\mathbf{P} = \mathbf 1$}

Let us show that the numerical scheme based on the equation
\begin{equation}\label{AppendixA1}
{}^{n+1} \mathbf C_{\text{i}} = {}^{n} \mathbf C_{\text{i}} +  2 \frac{\displaystyle
\xi }{\displaystyle \mathfrak{F}_2}
 \big( \mu \ {}^{n+1} \overline{\mathbf C}
 -  \frac{c_1}{2} \ {}^{n+1} \mathbf C_{\text{i}} \ {}^{n+1}\mathbf C_{\text{1i}}^{-1} \ {}^{n+1} \mathbf C_{\text{i}}
 -  \frac{c_2}{2} \ {}^{n+1} \mathbf C_{\text{i}} \ {}^{n+1}\mathbf C_{\text{2i}}^{-1} \
 {}^{n+1} \mathbf C_{\text{i}}  \big) + \beta \Delta t \ {}^{n+1} \mathbf C_{\text{i}} + \varepsilon \mathbf{1}
\end{equation}
is not weakly invariant. Toward that end, consider a reference change induced by $\textbf{F}_0$,
where $\text{det}(\textbf{F}_0) = 1$. Due to the reference change, the current right Cauchy-Green tensor ${}^{n+1} \textbf{C}$
transforms as follows
\begin{equation}\label{AppendixA2}
{}^{n+1} \overline{\mathbf C}^{\text{new}} = \textbf{F}^{-\text{T}}_0 \ {}^{n+1} \overline{\mathbf C} \ \textbf{F}^{-1}_0.
\end{equation}
The weak invariance can only take place if the internal variables are transformed according to
\begin{equation}\label{AppendixA3}
{}^{n} \textbf{C}_{\text{i}}^{\text{new}} = \textbf{F}^{-\text{T}}_0 \ {}^{n} \textbf{C}_{\text{i}} \ \textbf{F}^{-1}_0,
\end{equation}
\begin{equation}\label{AppendixA4}
{}^{n+1} \textbf{C}_{\text{1i}}^{\text{new}} = \textbf{F}^{-\text{T}}_0 \ {}^{n+1} \textbf{C}_{\text{1i}} \ \textbf{F}^{-1}_0, \quad
{}^{n+1} \textbf{C}_{\text{2i}}^{\text{new}} = \textbf{F}^{-\text{T}}_0 \ {}^{n+1} \textbf{C}_{\text{2i}} \ \textbf{F}^{-1}_0.
\end{equation}
Substituting these relations into the original equation \eqref{AppendixA1}, which is formulated with respect to unknown ${}^{n+1} \textbf{C}_{\text{i}}$, we arrive at
\begin{multline}\label{AppendixA5}
{}^{n+1} \mathbf C_{\text{i}} = \textbf{F}^{\text{T}}_0 \ {}^{n} \mathbf C_{\text{i}}^{\text{new}} \ \textbf{F}_0 +  2 \frac{\displaystyle
\xi }{\displaystyle \mathfrak{F}_2}
 \big( \mu \ \textbf{F}^{\text{T}}_0 {}^{n+1} \overline{\mathbf C}^{\text{new}} \ \textbf{F}_0
 -  \frac{c_1}{2} \ {}^{n+1} \mathbf C_{\text{i}} \ \textbf{F}^{-1}_0 ({}^{n+1} {\mathbf C}^{\text{new}}_{\text{1i}})^{-1} \
 \textbf{F}^{-\text{T}}_0 \ {}^{n+1} \mathbf C_{\text{i}} - \\
 -  \frac{c_2}{2} \ {}^{n+1} \mathbf C_{\text{i}} \ \textbf{F}^{-1}_0 ({}^{n+1} {\mathbf C}^{\text{new}}_{\text{2i}})^{-1} \
 \textbf{F}^{-\text{T}}_0 \ {}^{n+1} \mathbf C_{\text{i}}  \big) + \beta \Delta t \ {}^{n+1} \mathbf C_{\text{i}} + \varepsilon \mathbf{1}.
\end{multline}
Multiplying both sides of \eqref{AppendixA5} with $\textbf{F}^{-\text{T}}_0$ from the left and $\textbf{F}^{-1}_0$ from the right, and
introducing $\mathbf Z := \textbf{F}^{-\text{T}}_0 \ {}^{n+1} \mathbf C_{\text{i}} \ \textbf{F}^{-1}_0$, we obtain
\begin{equation}\label{AppendixA6}
\mathbf Z = {}^{n} \mathbf C_{\text{i}}^{\text{new}} +  2 \frac{\displaystyle
\xi }{\displaystyle \mathfrak{F}_2}
 \big( \mu \ {}^{n+1} \overline{\mathbf C}^{\text{new}}
 -  \frac{c_1}{2} \ \mathbf Z \ ({}^{n+1}\mathbf C_{\text{1i}}^{\text{new}})^{-1} \ \mathbf Z
 -  \frac{c_2}{2} \ \mathbf Z \ ({}^{n+1}\mathbf C_{\text{2i}}^{\text{new}})^{-1} \ \mathbf Z  \big) + \beta \Delta t \ \mathbf Z + \varepsilon \textbf{F}^{-\text{T}}_0 \textbf{F}^{-1}_0.
\end{equation}
Next, let ${}^{n+1} \textbf{C}_{\text{i}}^{\text{new}}$ be the solution of the new equation, which is formally obtained from equation
\eqref{AppendixA1} by replacing the old quantities by their new counterparts. This equation has almost the same from as \eqref{AppendixA6}. The
only difference lies in the last term on the right-hand side: $\varepsilon \textbf{1} \neq \varepsilon \textbf{F}^{-\text{T}}_0 \textbf{F}^{-1}_0$.
Therefore, these two equations have two different solutions: ${}^{n+1} \textbf{C}_{\text{i}}^{\text{new}} \neq \textbf{Z}$. Recalling the definition
of $\textbf{Z}$, we finally have
\begin{equation}\label{AppendixA7}
{}^{n+1} \textbf{C}_{\text{i}}^{\text{new}} \neq \textbf{F}^{-\text{T}}_0 \ {}^{n+1} \mathbf C_{\text{i}} \ \textbf{F}^{-1}_0.
\end{equation}
Thus, the weak invariance is violated. Numerical schemes of type \eqref{AppendixA1} should not be employed for discretization
of weakly invariant constitutive equations.

\subsection*{Weak invariance for $\mathbf{P} = {}^{n+1} \mathbf C_{\text{i}}$}

Now we prove the weak invariance property for the numerical scheme
\begin{multline}\label{AppendixA8}
{}^{n+1} \mathbf C_{\text{i}} = {}^{n} \mathbf C_{\text{i}} +  2 \frac{\displaystyle
\xi }{\displaystyle \mathfrak{F}_2}
 \big( \mu \ {}^{n+1} \overline{\mathbf C}
 -  \frac{c_1}{2} \ {}^{n+1} \mathbf C_{\text{i}} \ {}^{n+1}\mathbf C_{\text{1i}}^{-1} \ {}^{n+1} \mathbf C_{\text{i}}
 -  \frac{c_2}{2} \ {}^{n+1} \mathbf C_{\text{i}} \ {}^{n+1}\mathbf C_{\text{2i}}^{-1} \
 {}^{n+1} \mathbf C_{\text{i}}  \big) + \\ + \beta \Delta t \ {}^{n+1} \mathbf C_{\text{i}} + \varepsilon \ {}^{n+1} \mathbf C_{\text{i}},
\end{multline}
which is used to determine ${}^{n+1} \mathbf C_{\text{i}}$.
In order to prove this invariance, we need to assume
that all the input variables are transformed properly. More precisely, we suppose that
\eqref{AppendixA3} and \eqref{AppendixA4} hold true.
Substituting these relations into \eqref{AppendixA8}, we arrive at
\begin{multline}\label{AppendixA9}
{}^{n+1} \mathbf C_{\text{i}} = \textbf{F}^{\text{T}}_0 \ {}^{n} \mathbf C_{\text{i}}^{\text{new}} \ \textbf{F}_0 +  2 \frac{\displaystyle
\xi }{\displaystyle \mathfrak{F}_2}
 \big( \mu \ \textbf{F}^{\text{T}}_0 {}^{n+1} \overline{\mathbf C}^{\text{new}} \ \textbf{F}_0
 -  \frac{c_1}{2} \ {}^{n+1} \mathbf C_{\text{i}} \ \textbf{F}^{-1}_0 ({}^{n+1} {\mathbf C}^{\text{new}}_{\text{1i}})^{-1} \
 \textbf{F}^{-\text{T}}_0 \ {}^{n+1} \mathbf C_{\text{i}} - \\
 -  \frac{c_2}{2} \ {}^{n+1} \mathbf C_{\text{i}} \ \textbf{F}^{-1}_0 ({}^{n+1} {\mathbf C}^{\text{new}}_{\text{2i}})^{-1} \
 \textbf{F}^{-\text{T}}_0 \ {}^{n+1} \mathbf C_{\text{i}}  \big) + \beta \Delta t \ {}^{n+1} \mathbf C_{\text{i}} + \varepsilon {}^{n+1} \mathbf C_{\text{i}}.
\end{multline}
Multiplying both sides of \eqref{AppendixA9} with $\textbf{F}^{-\text{T}}_0$ from the left and $\textbf{F}^{-1}_0$ from the right,
we obtain the following equation for $\mathbf Z := \textbf{F}^{-\text{T}}_0 \ {}^{n+1} \mathbf C_{\text{i}} \ \textbf{F}^{-1}_0$
\begin{equation}\label{AppendixA10}
\mathbf Z = {}^{n} \mathbf C_{\text{i}}^{\text{new}} +  2 \frac{\displaystyle
\xi }{\displaystyle \mathfrak{F}_2}
 \big( \mu \ {}^{n+1} \overline{\mathbf C}^{\text{new}}
 -  \frac{c_1}{2} \ \mathbf Z \ ({}^{n+1}\mathbf C_{\text{1i}}^{\text{new}})^{-1} \ \mathbf Z
 -  \frac{c_2}{2} \ \mathbf Z \ ({}^{n+1}\mathbf C_{\text{2i}}^{\text{new}})^{-1} \ \mathbf Z  \big) + \beta \Delta t \ \mathbf Z + \varepsilon \mathbf Z.
\end{equation}
This equation exactly coincides with \eqref{AppendixA8}, where the old quantities are replaced by the new ones. Thus, ${}^{n+1} \textbf{C}_{\text{i}}^{\text{new}} = \textbf{Z}$
\begin{equation}\label{AppendixA10}
{}^{n+1} \textbf{C}_{\text{i}}^{\text{new}} = \textbf{F}^{-\text{T}}_0 \ {}^{n+1} \mathbf C_{\text{i}} \ \textbf{F}^{-1}_0,
\end{equation}
which is exactly the required weak invariance.

\section*{Appendix B. Analysis of round-off errors}

Let us estimate the round-off errors which appear during the step-by-step evaluation of \eqref{ReducedEvolut5}.
The square root in \eqref{ReducedEvolut5} can be computed exactly
up to machine precision $\epsilon$. For simplicity assume that this is the only source of errors.
Thus, the step-by-step evaluation of \eqref{ReducedEvolut5} yields
\begin{equation}\label{Contamination}
\mathbf Y \leftarrow \frac{\displaystyle
 \mathfrak{F}_2}{\displaystyle 2 \xi c} \Big[ -z \textbf{1} + \big(z^2 \textbf{1} +  4 \frac{\displaystyle
\xi c}{\displaystyle \mathfrak{F}_2} \textbf{A} \big)^{1/2} + \boldsymbol{\epsilon} \Big], \quad \text{where} \quad  \| \boldsymbol{\epsilon}  \| \approx \epsilon.
\end{equation}
In the machine arithmetics we obtain
\begin{equation}\label{Contamination0}
\mathbf Y^{\text{machine}} \leftarrow \mathbf Y^{\text{exact}} + \frac{\displaystyle
 \mathfrak{F}_2}{\displaystyle 2 \xi c}  \boldsymbol{\epsilon}.
\end{equation}
In other words, $\mathbf Y^{\text{machine}} $ is contaminated by the error $\frac{\displaystyle
 \mathfrak{F}_2}{\displaystyle 2 \xi c}  \boldsymbol{\epsilon}$, which
can become very large for
small $2 \frac{ \displaystyle \xi c}{\displaystyle \mathfrak{F}_2}$.

Now let us analyze the alternative equation \eqref{ReducedEvolut5notCont}.
Its step-by-step evaluation is equivalent to the following chain
\begin{equation}\label{Chain0}
\mathbf{Y}_0 \leftarrow \big(z^2 \textbf{1} +  4 \frac{\displaystyle
\xi c}{\displaystyle \mathfrak{F}_2} \textbf{A} \big)^{1/2}, \ \mathbf{Y}_1 \leftarrow \mathbf{Y}_0 + z \mathbf{1}, \
\mathbf{Y}_2 \leftarrow \mathbf{Y}^{-1}_1, \ \mathbf{Y} \leftarrow 2 \mathbf{A} \mathbf{Y}_2.
\end{equation}
Again, the square root is computed with a small error $\boldsymbol{\epsilon}$.
Thus, we obtain
\begin{equation}\label{Contamination2}
\mathbf{Y}^{\text{machine}}_0 \leftarrow \mathbf{Y}^{\text{exact}}_0 + \boldsymbol{\epsilon}, \
\mathbf{Y}^{\text{machine}}_1 \leftarrow \mathbf{Y}^{\text{exact}}_1 + \boldsymbol{\epsilon}.
\end{equation}
Using the Neumann series, it can be easily shown that
\begin{equation}\label{Contamination3}
(\mathbf{Y}^{\text{exact}}_1 + \boldsymbol{\epsilon})^{-1} =
(\mathbf{Y}^{\text{exact}}_1)^{-1} - (\mathbf{Y}^{\text{exact}}_1)^{-1} \ \boldsymbol{\epsilon} \ (\mathbf{Y}^{\text{exact}}_1)^{-1} +
O(\epsilon^2).
\end{equation}
Thus, neglecting the higher-order terms, we have
\begin{equation}\label{Contamination4}
\mathbf{Y}^{\text{machine}}_2 \leftarrow \mathbf{Y}^{\text{exact}}_2 - \mathbf{Y}^{\text{exact}}_2 \ \boldsymbol{\epsilon}  \ \mathbf{Y}^{\text{exact}}_2.
\end{equation}
Finally,
\begin{equation}\label{Contamination4}
\mathbf{Y}^{\text{machine}} \leftarrow \mathbf{Y}^{\text{exact}} - 2 \mathbf{A}
\mathbf{Y}^{\text{exact}}_2 \ \boldsymbol{\epsilon}  \ \mathbf{Y}^{\text{exact}}_2 =
\mathbf{Y}^{\text{exact}} - \mathbf{Y}^{\text{exact}} \ \boldsymbol{\epsilon}  \ \mathbf{Y}^{\text{exact}}_2.
\end{equation}
Since $\mathbf{Y}^{\text{exact}}_2$ is bounded, the small error $\boldsymbol{\epsilon}$ \emph{is not multiplied} by a large factor as it was the
case for \eqref{ReducedEvolut5}.
In other words, relation \eqref{ReducedEvolut5notCont} can be evaluated step-by-step.

\section*{Appendix C. Estimation of $z$}

Let us derive the estimation \eqref{EstimOfz} for the unknown parameter $z$.
First, we recall equation \eqref{QuadrEquation} and the incompressibility relation formulated in terms of $\mathbf Y$
\begin{equation}\label{Appendix1}
z \ \mathbf Y = \mathbf A -  \frac{\displaystyle
\xi c}{\displaystyle \mathfrak{F}_2} \mathbf Y^2, \quad \det (\textbf{Y}) = \det (\mathbf{\Phi}).
\end{equation}
In order to rewrite the first relation in a more compact form, we introduce the abbreviation $\xi':=\frac{\displaystyle c \xi}{\displaystyle \mathfrak{F}_2}$
\begin{equation}\label{Appendix2}
z \ \mathbf Y = \mathbf A -  \xi' \mathbf Y^2.
\end{equation}
This equation yields $\mathbf Y$ as a function of $z$ and $\xi'$.
Let us consider its expansion in Taylor series for small  $\xi'$
\begin{equation}\label{Appendix3}
\mathbf Y = \tilde{\mathbf Y}(z,\xi') = \frac{1}{z} \mathbf A - \frac{\xi'}{z^3} {\mathbf A}^2 + O(\xi'^2), \quad
\mathbf Y|_{\xi' =0} = \frac{1}{z} \mathbf A.
\end{equation}
Since the parameter $z$ was introduced
to enforce the incompressibility, $z$ is estimated using the incompressibility relation $\eqref{Appendix1}_2$, which yields
\begin{equation}\label{Appendix4}
z = \tilde{z} (\xi'), \quad z_0 := \tilde{z} (0) = \Big(\frac{\det \mathbf A}{\det{\mathbf \Phi}}\Big)^{1/3}.
\end{equation}
Employing the implicit function theorem, we have
\begin{equation}\label{Appendix5}
\frac{d \tilde{z} (\xi')}{d \xi'}|_{\xi' =0} = - \frac{\partial \det \tilde{\mathbf Y}(z,\xi')}{\partial \xi'}|_{z=z_0, \xi' =0}
\ \Big(\frac{\partial \det \tilde{\mathbf Y}(z,\xi')}{\partial z}|_{z=z_0, \xi' =0}\Big)^{-1}.
\end{equation}
Next, adopting the Jacobi formula and differentiating $\eqref{Appendix3}_1$, we obtain
\begin{equation}\label{Appendix52}
\frac{\partial \det \tilde{\mathbf Y}(z,\xi')}{\partial \xi'}|_{z=z_0, \xi' =0} =
\det \tilde{\mathbf Y}(z_0,0) \ (\tilde{\mathbf Y}(z_0,0))^{-1} : \Big(-\frac{\mathbf A^2}{z_0^3}\Big),
\end{equation}
\begin{equation}\label{Appendix53}
\frac{\partial \det \tilde{\mathbf Y}(z,\xi')}{\partial z}|_{z=z_0, \xi' =0} =
\det \tilde{\mathbf Y}(z_0,0) \ (\tilde{\mathbf Y}(z_0,0))^{-1} : \Big(- \frac{\mathbf A}{z_0^2}\Big).
\end{equation}
Substituting this into \eqref{Appendix5}, after some algebraic computations we arrive at
\begin{equation}\label{Appendix6}
\frac{d \tilde{z} (\xi')}{d \xi'}|_{\xi' =0} = - \frac{\text{tr} \mathbf A}{3 z_0}, \quad
\tilde{z}(\xi') = z_0 - \frac{\text{tr} \mathbf A}{3 z_0} \xi' + O(\xi'^2).
\end{equation}
Finally, we have
\begin{equation}\label{Appendix7}
z = z_0 - \frac{\text{tr} \mathbf A}{3 z_0} \frac{c \ \xi}{\mathfrak{F}_2} + O\Big(\Big(\frac{c \ \xi}{\mathfrak{F}_2}\Big)^2\Big).
\end{equation}
Obviously, this estimation is exact for $c=0$ and finite $\xi$.

\section*{Appendix D. Weak invariance of a certain push-forward operation}

\begin{figure}\centering
\psfrag{A}[m][][1][0]{$\overline{{}^{n}\mathbf{C}}$}
\psfrag{B}[m][][1][0]{${}^{n}\mathbf{C}_{\text{i}}$}
\psfrag{C}[m][][1][0]{${}^{\text{est}}\mathbf{C}_{\text{i}}$}
\psfrag{D}[m][][1][0]{$\overline{{}^{n+1}\mathbf{C}}$}
\psfrag{E}[m][][1][0]{${}^{n}\mathbf{C}^{\text{new}}_{\text{i}}$}
\psfrag{F}[m][][1][0]{${}^{\text{est}}\mathbf{C}^{\text{new}}_{\text{i}}$}
\psfrag{G}[m][][1][0]{$\overline{{}^{n}\mathbf{C}^{\text{new}}}$}
\psfrag{H}[m][][1][0]{$\overline{{}^{n+1}\mathbf{C}^{\text{new}}}$}
\psfrag{K}[m][][1][0]{$\textbf{F}^{-\text{T}}_0 \ (\cdot) \ \textbf{F}^{-1}_0$}
\psfrag{L}[m][][1][0]{?}
\psfrag{M}[m][][1][0]{$\mathbf{F}^{-\text{T}}_{\text{sh}} \ (\cdot) \ \mathbf{F}^{-1}_{\text{sh}}$}
\psfrag{N}[m][][1][0]{$(\mathbf{F}^{\text{new}}_{\text{sh}})^{-\text{T}} \ (\cdot) \ (\mathbf{F}^{\text{new}}_{\text{sh}})^{-1}$}
\scalebox{0.95}{\includegraphics{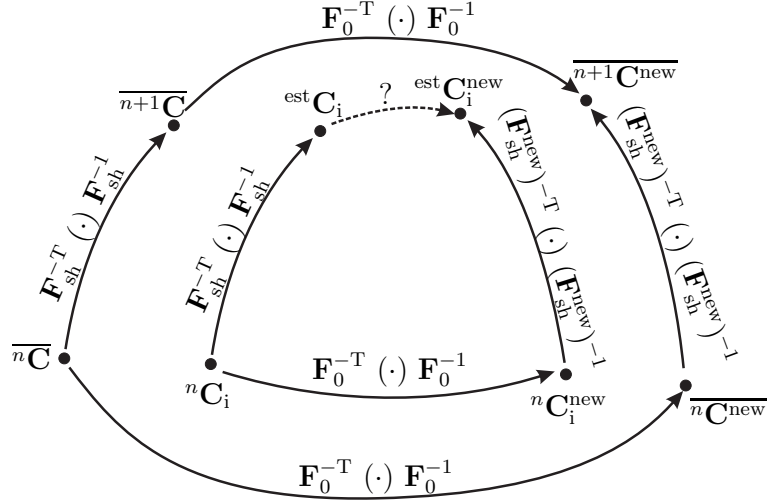}} \caption{Invariance of the push-forward operation
\eqref{DefineShift} used to estimate ${}^{n+1}\mathbf{C}_{\text{i}}$. Quantities on the left-hand side of the figure operate on the original
reference configuration; quantities on the right-hand side operate on the new reference.
\label{AppendixC}}
\end{figure}

Let us check that the push-forward operation \eqref{DefineShift} preserves the weak invariance of the solution.
Toward that end consider the quantities ${}^{n+1}\mathbf{C}$, ${}^{n}\mathbf{C}$, ${}^{n} \mathbf{C}_{\text{i}}$
which operate on the original reference configuration (see Figure \ref{AppendixC}).
Operation \eqref{DefineShift} can be rewritten in the following way
\begin{equation}\label{AppendixC1}
{}^{\text{est}}  \mathbf{C}_{\text{i}} := \mathbf{F}^{-\text{T}}_{\text{sh}} \ {}^{n}\mathbf{C}_{\text{i}} \ \mathbf{F}^{-1}_{\text{sh}} =
(\overline{{}^{n+1}  \mathbf{C} \  {}^{n}  \mathbf{C}^{-1}})^{1/2} \ {}^{n}\mathbf{C}_{\text{i}}
\ (\overline{{}^{n}  \mathbf{C}^{-1} \ {}^{n+1}  \mathbf{C}})^{1/2},
\end{equation}
where the push-forward $\mathbf{F}^{-\text{T}}_{\text{sh}} \ (\cdot) \ \mathbf{F}^{-1}_{\text{sh}}$
brings $\overline{{}^{n}\mathbf{C}}$ to $\overline{{}^{n+1}\mathbf{C}}$.
Let $\textbf{F}_0$ be the reference change, such that $\text{det}(\textbf{F}_0) = 1$. The quantities with respect to the new reference configuration are computed as follows
\begin{equation}\label{AppendixC2}
{}^{n+1}\mathbf{C}^{\text{new}} = \textbf{F}^{-\text{T}}_0 \ {}^{n+1}\mathbf{C} \ \textbf{F}^{-1}_0, \quad
{}^{n}\mathbf{C}^{\text{new}} = \textbf{F}^{-\text{T}}_0 \ {}^{n}\mathbf{C} \ \textbf{F}^{-1}_0, \quad
{}^{n}\mathbf{C}_{\text{i}}^{\text{new}} = \textbf{F}^{-\text{T}}_0 \ {}^{n}\mathbf{C}_{\text{i}} \ \textbf{F}^{-1}_0.
\end{equation}
Specifying \eqref{AppendixC1} for the quantities on the new reference, we have
\begin{equation}\label{AppendixC3}
{}^{\text{est}}  \mathbf{C}^{\text{new}}_{\text{i}} :=
(\mathbf{F}_{\text{sh}}^{\text{new}})^{-\text{T}} \ {}^{n}\mathbf{C}^{\text{new}}_{\text{i}} \ (\mathbf{F}^{\text{new}}_{\text{sh}})^{-1} =
(\overline{{}^{n+1}  \mathbf{C}^{\text{new}} \  ({}^{n}  \mathbf{C}^{\text{new}})^{-1}})^{1/2} \ {}^{n}\mathbf{C}^{\text{new}}_{\text{i}}
\ (\overline{ ({}^{n}\mathbf{C}^{\text{new}})^{-1} \ {}^{n+1}\mathbf{C}^{\text{new}}})^{1/2},
\end{equation}
where the push-forward $(\mathbf{F}_{\text{sh}}^{\text{new}})^{-\text{T}} \ (\cdot) \ (\mathbf{F}^{\text{new}}_{\text{sh}})^{-1}$
brings $\overline{{}^{n}\mathbf{C}^{\text{new}}}$ to $\overline{{}^{n+1}\mathbf{C}^{\text{new}}}$.
The situation is summarized in Figure \ref{AppendixC}.

Our goal now is to clarify the interrelationship between
${}^{\text{est}}\mathbf{C}^{\text{new}}_{\text{i}}$ and ${}^{\text{est}}\mathbf{C}_{\text{i}}$.
Substituting \eqref{AppendixC2} into \eqref{AppendixC3} and using the identities $(\textbf{F}^{-\text{T}}_0 \mathbf{A} \textbf{F}^{\text{T}}_0)^{1/2} =
\textbf{F}^{-\text{T}}_0 \ \mathbf{A}^{1/2} \ \textbf{F}^{\text{T}}_0$ as well as
$(\textbf{F}_0 \mathbf{A} \textbf{F}^{-1}_0)^{1/2} =
\textbf{F}_0 \ \mathbf{A}^{1/2} \ \textbf{F}^{-1}_0$, we arrive at
\begin{equation}\label{AppendixC4}
{}^{\text{est}}  \mathbf{C}^{\text{new}}_{\text{i}} = \textbf{F}^{-\text{T}}_0 \
(\overline{{}^{n+1}  \mathbf{C} \  {}^{n}  \mathbf{C}^{-1}})^{1/2} \ {}^{n}\mathbf{C}_{\text{i}}
\ (\overline{{}^{n}  \mathbf{C}^{-1} \ {}^{n+1}  \mathbf{C}})^{1/2} \ \textbf{F}^{-1}_0 = \textbf{F}^{-\text{T}}_0
\  {}^{\text{est}}  \mathbf{C}_{\text{i}}  \  \textbf{F}^{-1}_0.
\end{equation}
This relation proves the weak invariance of push-forward \eqref{AppendixC1}.

Alternatively, one may consider another push-forward operation, defined thorough
\begin{equation}\label{Appendix534}
{}^{\text{est}}  \mathbf{C}_{\text{i}} \gets \tilde{\mathbf{F}}^{-\text{T}}_{\text{sh}} \ {}^{n}  \mathbf{C}_{\text{i}} \ \tilde{\mathbf{F}}^{-1}_{\text{sh}}, \ \text{where} \
 \tilde{\mathbf{F}}_{\text{sh}} := (\overline{{}^{n+1}  \mathbf{C}})^{-1/2} \  (\overline{{}^{n}  \mathbf{C}})^{1/2}.
\end{equation}
Just as \eqref{DefineShift}, it brings $\overline{{}^{n}\mathbf{C}}$ to $\overline{{}^{n+1}\mathbf{C}}$, but
such a transformation is not invariant under the change of the reference configuration:
\begin{equation}\label{Appendix6}
{}^{\text{est}}  \mathbf{C}^{\text{new}}_{\text{i}} \neq \textbf{F}^{-\text{T}}_0
\  {}^{\text{est}}  \mathbf{C}_{\text{i}}  \  \textbf{F}^{-1}_0.
\end{equation}
Thus, operations of type \eqref{Appendix534} should not be used in numerical computations.

\section*{Appendix E. EBMSC and EM}

\subsubsection*{Discretization of differential equations using EBMSC and EM}

Let us consider the initial value problem for a system of
nonlinear ordinary differential equations
\begin{equation*}\label{difur}
\dot{\mathbf A} (t) = \mathbf{f} (\mathbf A(t), t) \mathbf A (t), \quad
\mathbf A(0)=\mathbf A^0, \quad \det(\mathbf A^0)=1.
\end{equation*}

By ${}^n \mathbf A, {}^{n+1} \mathbf A$ denote numerical solutions respectively at $t_n$ and $t_{n+1}$.
The classical Euler Backward method (EBM) uses the equation with respect to the unknown ${}^{n+1} \mathbf A^{\text{EBM}}$ \cite{SimoHughes, SimMieh}:
\begin{equation}\label{Eulcl}
{}^{n+1} \mathbf A^{\text{EBM}} = \big[ \mathbf 1 - \Delta t \ \mathbf{f} ({}^{n+1} \mathbf A^{\text{EBM}}, t_{n+1}) \big]^{-1}
\ {}^n \mathbf A.
\end{equation}
If the incompressibility restriction $\det ({}^{n+1} \mathbf A) \equiv 1$ needs to be enforced, a modified method, called Euler Backward method with subsequent correction of incompressibility (EBMSC), can be considered \cite{ShutovLandgraf2013}:
\begin{equation}\label{EBMSC}
{}^{n+1} \mathbf A^{\text{EBMSC}} := \overline{{}^{n+1} \mathbf A^{\text{EBM}}}.
\end{equation}
The exponential method (EM) is based on the equation with respect to ${}^{n+1} \mathbf A^{\text{EM}}$ \cite{WebAnan, MiStei, Simo}:
\begin{equation}\label{Expo}
{}^{n+1} \mathbf A^{\text{EM}} = \exp\big(\displaystyle \Delta t \ \mathbf{f} ({}^{n+1} \mathbf A^{\text{EM}}, t_{n+1})\big) \ {}^n \mathbf A.
\end{equation}
In contrast to the classical EBM, EM exactly preserves the incompressibility.
Under some general assumptions, both methods preserve the symmetry of $\mathbf A$ (see \cite{ShutovKrVisc}).

\subsubsection*{Application to the model of Shutov and Krei\ss ig}

For the material model under consideration, two conventional numerical procedures can be obtained using EBMSC and EM in the following way.
First, we discretize the evolution equations \eqref{prob1} -- \eqref{prob22} using EBMSC or EM.
Next, the scalar evolution equations \eqref{prob3} are discretized using the Euler Backward method as discussed in Section 3.1.
Within a nested procedure, the resulting system of algebraic equations is solved numerically for a fixed
$\xi := \Delta t \ {}^{n+1} \lambda_{\text{i}}$ with respect to unknown
$\mathbf C_{\text{i}}$, $\mathbf C_{\text{1i}}$, $\mathbf C_{\text{2i}}$, $s$, and $s_{\text{d}}$.
Implementing these functional dependencies, the
unknown inelastic strain increment $\xi$ is found using the predictor-corrector scheme, by
solving the consistency condition
\begin{equation}\label{AppendixE}
\xi \eta = \displaystyle \Delta t \Big\langle \frac{\displaystyle 1}{\displaystyle f_0}
f \Big\rangle^{m},
\end{equation}
where $f$ is a function of $\xi$. For more details, the reader is referred to \cite{ShutovKrVisc, ShutovKrKoo}.

\end{document}